\documentclass[11pt,american,english]{article}
\usepackage[T1]{fontenc}
\usepackage[latin9]{inputenc}
\usepackage[a4paper]{geometry}
\usepackage{babel}
\usepackage{textcomp}
\usepackage{amsthm}
\usepackage{amsmath}
\usepackage{amssymb}
\usepackage{setspace}
\usepackage{esint}
\usepackage[unicode=true,pdfusetitle,
 bookmarks=true,bookmarksnumbered=true,bookmarksopen=true,bookmarksopenlevel=3,
 breaklinks=false,pdfborder={0 0 1},backref=false,colorlinks=false]
 {hyperref}
\usepackage{breakurl}

\makeatletter
\theoremstyle{plain}
\newtheorem{thm}{\protect\theoremname}
  \theoremstyle{remark}
  \newtheorem{rem}[thm]{\protect\remarkname}
  \theoremstyle{definition}
  \newtheorem{example}[thm]{\protect\examplename}
  \theoremstyle{definition}
  \newtheorem{xca}[thm]{\protect\exercisename}

\usepackage{multirow}

\usepackage{arydshln}

\usepackage{lmodern}
\newcommand{\1}{\mbox{1\hspace{-1mm}I}}
\numberwithin{equation}{section}

\usepackage{calc}
\makeatother

  \addto\captionsamerican{\renewcommand{\examplename}{Example}}
  \addto\captionsamerican{\renewcommand{\exercisename}{Exercise}}
  \addto\captionsamerican{\renewcommand{\remarkname}{Remark}}
  \addto\captionsamerican{\renewcommand{\theoremname}{Theorem}}
  \addto\captionsenglish{\renewcommand{\examplename}{Example}}
  \addto\captionsenglish{\renewcommand{\exercisename}{Exercise}}
  \addto\captionsenglish{\renewcommand{\remarkname}{Remark}}
  \addto\captionsenglish{\renewcommand{\theoremname}{Theorem}}
  \providecommand{\examplename}{Example}
  \providecommand{\exercisename}{Exercise}
  \providecommand{\remarkname}{Remark}
\providecommand{\theoremname}{Theorem}

\begin{document}
\selectlanguage{american}%
\global\long\def\1{\mbox{1\hspace{-1mm}I}}

\selectlanguage{english}%

\title{ON THE INTERPRETATION OF THE MASTER EQUATION}

\author{Alain Bensoussan
\thanks{axb046100@utdallas.edu}\\
International Center for Decision and Risk Analysis\\
Jindal School of Management, University of Texas at Dallas%
\\
Jens Frehse 
\thanks{mathfrehse@googlemail.com}\\
Institute for Applied Mathematics,University of Bonn
\\
Sheung Chi Phillip Yam
\thanks{scpyam@sta.cuhk.edu.hk}\\
Department of Statistics, The Chinese University of Hong Kong
}

\maketitle
\footnotetext{The first author-Alain Bensoussan is also affiliated with the College of Science and Engineering, Systems Engineering and Engineering Management, City University Hong Kong. Research supported by the Hong Kong RGC GRF 500113, the National Science Foundation under grant DMS-1303775, and the Research Grants Council of the Hong Kong Special Administrative Region(CityU 500113). The third author-Phillip Yam acknowledges the financial supports from The Hong Kong RGC GRF 404012 with the project title: Advanced Topics In Multivariate Risk Management In Finance And Insurance, and Direct Grant for Research 2014/15 with project code: 4053141 offered by CUHK.
}
\begin{abstract}
Since its introduction by P.L. Lions in his lectures and seminars
at the {\it{College de France}}, see \cite{PLL}, and also the very helpful
notes of Cardialaguet \cite{PCA} on Lions' lectures, the Master Equation
has attracted a lot of interest, and various points of view have been
expressed, see for example Carmona-Delarue \cite{CAD}, Bensoussan-Frehse-Yam
\cite{BFY1}, Buckdahn-Li-Peng-Rainer \cite{BLPR}. There are several
ways to introduce this type of equation; and in those mentioned works, they involve an argument which
is a probability measure, while P.L. Lions has recently proposed the idea of
working with the Hilbert space of square integrable random variables. Hence writing the equation is an issue; while another issue is
its origin. In this article, we discuss all these various aspects, and our modeling argument relies heavily on a seminar at {\it{College de France}} delivered
by P.L. Lions on November 14, 2014. 
\end{abstract}

\section{INTRODUCTION }

Our objective in this paper is to present and discuss some of the
main ideas and concepts related to the Master Equation, first introduced
by P.L. Lions. Since the major source of information is only available
on video, our present presentation has also the objective to be a survey.
The notes of Cardialaguet on Lions's lectures are extremely valuable; however, they are not in the form as a survey. Besides, we cover some additional
material, which is obtained more recently, including the seminar talk delivered by P.L. Lions
on November 14, 2014, and our own contribution. The very exciting
results of existence and uniqueness of the solution of the Master
Equation mentioned by P.L. Lions in this seminar are unfortunately
not covered here, since only the outline of their proofs were sketched in his talk. Instead, our present motivation is to provide a guidebook. A situation, which looks clear, is when one writes down
a Bellman equation; it is the case for mean field type control
problems, or mean field games which are equivalent to mean field type
control problems. The Master Equation is then obtained by differentiating
a Bellman equation, see \cite{BFY1}. One can have several writings but they all correspond
to the same solution. Nevertheless, there are other cases, in which the Master
Equation is not inherited from a Bellman equation. One important motivation
in such cases is the fact that the Master Equation is a way to uncouple
a forward-backward system of deterministic or stochastic P.D.E.s corresponding
to mean field games or mean field type control problems. Another point of view
is the fact that mean field games approximate differential games
with a large number of players. The connection between the Master
Equation and the large stochastic differential game also poses an interesting question.
Finally, the linear quadratic model is quite useful, since one can often obtain explicit
formulae, and so facilitates different comparisons.

\section{FUNCTIONALS OF PROBABILITY MEASURES}

\subsection{DERIVATIVES}

In all problems connected with mean field theory, one naturally works with functionals of probability measures
on $\mathbb{R}^{n}$. An important issue is to define a convenient functional
space. The natural space is the Wasserstein space $\mathcal{P}^{2}$
. But in order to differentiate, this structure is not convenient since it is generally not a vector space.
Two approaches are possible. Let $m$ be a probability measure on
$\mathbb{R}^{n}.$ If it has a density, we shall use the same notation for
the density, except that $m\equiv m(x),$ with $x\in \mathbb{R}^{n}$. We
then assume that $m\in L^{2}(\mathbb{R}^{n}).$ If we consider a functional
$F(m),$ we assume that it is defined on $\mathcal{P}^{2}$ and at
the same time on $L^{2}(\mathbb{R}^{n}).$ We can use the concept of Gâteaux
differentiability on $L^{2}(\mathbb{R}^{n}).$ If $F:L^{2}(\mathbb{R}^{n})\rightarrow \mathbb{R},$
then

\[
\frac{d}{d\theta}F\left(m+\theta\tilde{m}\right)\bigg|_{\theta=0}=\int_{\mathbb{R}^{n}}\frac{\partial F}{\partial m}(m)(\xi)\tilde{m}(\xi)d\xi,
\]
with $\xi\mapsto\frac{\partial F}{\partial m}(m)(\xi)$ being in
$L^{2}(\mathbb{R}^{n}).$ The second derivative $\frac{\partial^{2}F}{\partial m^{2}}\left(m\right)$ is a linear map from $L^{2}(\mathbb{R}^{n})$ into $L^{2}(\mathbb{R}^{n})^{\otimes 2}$ such that 

\[
\frac{d}{d\theta}\frac{\partial F}{\partial m}(m+\theta\tilde{m})(\xi)\bigg|_{\theta=0}=\int_{\mathbb{R}^{n}}\frac{\partial^{2}F}{\partial m^{2}}(m)(\xi,\eta)\tilde{m}(\eta)d\eta,
\]

where the function $\frac{\partial^{2}F}{\partial m^{2}}(m)(\xi,\eta)$ is symmetric in $\left(\xi,\eta\right)$.

We can state the second order Taylor's formula 

\begin{equation}
F(m+\tilde{m})=F(m)+\int_{\mathbb{R}^{n}}\frac{\partial F}{\partial m}(m)(\xi)\tilde{m}(\xi)d\xi+\int_{0}^{1}\int_{0}^{1}\lambda\int_{\mathbb{R}^{n}}\int_{\mathbb{R}^{n}}\frac{\partial^{2}F}{\partial m^{2}}(m+\lambda\mu\tilde{m})(\xi,\eta)\tilde{m}(\xi)\tilde{m}(\eta)d\xi d\eta d\lambda d\mu\label{eq:1-1}
\end{equation}

The second possibility, introduced by P.L. Lions is to consider a
probability space $\Omega,\mathcal{A},P$ on which one can construct
a random variable $X$ on $\mathbb{R}^{n}$, such that the probability law
of $X,$ denoted by $\mathcal{L}_{X}=m.$ We then write $F(X)$ and
consider that $F$ is defined on the Hilbert space $\mathcal{H}=L^{2}(\Omega,\mathcal{A},P;\mathbb{R}^{n}).$
We use the same notation $F,$ to save notation. The fact that $F(X)=F(m)=F(\mathcal{L}_{X})$
means that the functional $F(X)$ depends of $X$ only through its
probability. This implies specific aspects, as we shall see. An obvious
advantage of this approach compared to the previous one is that it
is more general, since the probability measure does not need to have
a density. Besides, because we use random variables, we keep the fact
that $m$ is a probability, and is not a general function in $L^{2}(\mathbb{R}^{n}).$
We can then use the concept of Gâteaux differentiability on $\mathcal{H}$
. We denote the scalar product in $\mathcal{H}$ by $(X,Y).$ Of course
$(X,Y)=\mathbb{E}X.Y,$ but it is good to consider $\mathcal{H}$ as a general
Hilbert space. We then have 

\begin{equation}
\frac{d}{d\theta}F(X+\theta Y)|_{\theta=0}=(DF(X),Y)\label{eq:1-11-1}
\end{equation}
 and $DF(X)$$\in\mathcal{H}.$ We then proceed to define the 2nd
derivative 

\[
\frac{d}{d\theta}(DF(X+\theta Y),Z)|_{\theta=0}=(D^{2}F(X)Y,Z)=D^{2}F(X)(Y,Z)
\]

in which $D^{2}F(X)$$\in\mathcal{L}(\mathcal{H};\mathcal{H})$ is
identified to a symmetric bilinear form on $\mathcal{H}$. We can
then write as for (\ref{eq:1-1}) 

\begin{equation}
F(X+Y)=F(X)+(DF(X),Y)+\int_{0}^{1}\int_{0}^{1}\lambda D^{2}F(X+\lambda\mu Y)(Y,Y)\label{eq:1-2}
\end{equation}

The first question is what is the relation between $DF(X)$ and $\dfrac{\partial F(m)}{\partial m}(x)$
?

To answer this question we consider a random variable $Y$ and assume
that the pair $X,Y$ has a joint probability density $\pi(x,y).$
By definition the marginal of $X$ has the density $m(x).$ We define
the stochastic process $\xi(t)=X+tY.$ It is easy to check that for
$t>0$ , $\xi(t)$ has a probability density $m(x,t)$ given by the
formula 

\[
m(x,t)=\int\pi(x-ty,y)dy
\]
 so 

\[
\frac{\partial m}{\partial t}(x,t)=-\text{div}_{x}\int\pi(x-ty,y)\, ydy
\]

But then

\begin{align*}
\frac{d}{dt}F(m(t)) & =\int\frac{\partial F}{\partial m}(m(t))(x)\frac{\partial m}{\partial t}(x,t)dx\\
 & =-\int\frac{\partial F}{\partial m}(m(t))(x)\text{div}_{x}\int\pi(x-ty,y)\, ydydx
\end{align*}
 therefore 

\begin{align*}
\frac{d}{dt}F(m(t))|_{t=0} & =-\int\frac{\partial F}{\partial m}(m)(x)\text{div}_{x}\int\pi(x,y)\, ydydx\\
 & =\int D_{x}\frac{\partial F}{\partial m}(m)(x).(\int\pi(x,y)\, ydy)dx
\end{align*}
 From the definition of $\pi(x,y)$ this can be written as 

\[
\frac{d}{dt}F(m(t))|_{t=0}=\mathbb{E}[D_{x}\frac{\partial F}{\partial m}(m)(X).Y]
\]
On the other hand, since $F(m(t))=F(\xi(t))$ we have also 

\[
\frac{d}{dt}F(m(t))|_{t=0}=\frac{d}{dt}F(\xi(t))|_{t=0}=\mathbb{E}[DF(X).Y]
\]
 and since $Y$ is arbitrary, we have obtained the formula 

\begin{equation}
DF(X)=D_{x}\frac{\partial F}{\partial m}(m)(X)\label{eq:1-3}
\end{equation}

This proves a result which can be obtained in a more general set-up.
The derivative $DF(X)$ can be computed through a function of $m$
and $x$, denoted (using the notation of \cite{BLPR})$\partial_{m}F(m)(x)$,
with the formula 

\begin{equation}
DF(X)=\partial_{m}F(m)(X)\label{eq:1-4}
\end{equation}

Of course from the previous calculation, we can state that when the
probability has a density denoted by $m,$ then 

\begin{equation}
\partial_{m}F(m)(x)=D_{x}\frac{\partial F}{\partial m}(m)(x)\label{eq:1-5}
\end{equation}

We can apply this result to the functional $\tilde{F}(X)=\mathbb{E}F(X,m),$
which means 

\begin{equation}
\tilde{F}(X)=\int_{\mathbb{R}^{n}}F(x,m)m(x)dx\label{eq:1-31}
\end{equation}

We see easily , by formula (\ref{eq:1-3}) that 

\begin{equation}
D\tilde{F}(X)=D_{x}F(X,m)+D_{x}\int_{\mathbb{R}^{n}}\frac{\partial}{\partial m}F(\xi,m)(X)m(\xi)d\xi\label{eq:1-32}
\end{equation}

We will have a more complex situation to handle. We consider a function
$V(x,m):$$\mathbb{R}^{n}\times L^{2}(\mathbb{R}^{n})\rightarrow \mathbb{R}^{n}$. Next we set
$\mathcal{U}(X)=V(X,\mathcal{L}_{X})$ . Then $\mathcal{U}$ maps
$\mathcal{H}$ into itself. 

We want to compute its gradient $D\mathcal{U}(X)$ which is a map
from $\mathcal{H}$ into $\mathcal{L}(\mathcal{H};\mathcal{H}).$
If $Z\in\mathcal{H}$ then for any $W$ in $\mathcal{H}$ we have 

\begin{align}
(D\mathcal{U}(X)Z,W) & =(D_{X}(\mathcal{U}(X),W),Z)=\label{eq:1-34}
\end{align}

\[
=\mathbb{E}\sum_{j}Z_{j}\frac{\partial V_{i}}{\partial x_{j}}(X,\mathcal{L}_{X})W_{i}+(D_{X}\mathbb{E}\sum_{i}V_{i}(Y,\mathcal{L}_{X})\tilde{W}_{i},Z)
\]

in which $Y,\tilde{W}$ is a copy of $X,W$ independent of $X,Z.$
This implies that the dependence in $X$ in the second term is only
trough $\mathcal{L}_{X}.$ We then have 

\[
(D_{X}\mathbb{E}\sum_{i}V_{i}(Y,\mathcal{L}_{X})\tilde{W}_{i},Z)=\sum_{j}\mathbb{E}_{XZ}\, Z_{j}\frac{\partial}{\partial x_{j}}\frac{\partial}{\partial m}[\, \mathbb{E}_{Y\tilde{W}}\,\sum_{i}V_{i}(Y,\mathcal{L}_{X})(X)\tilde{W}_{i}]
\]
 with the notation $\mathbb{E}_{Y\tilde{W}}\,\sum_{i}V_{i}(Y,\mathcal{L}_{X})(X)\tilde{W}_{i}=\mathbb{E}\,$$\sum_{i}V_{i}(Y,\mathcal{L}_{X})(x)\tilde{W}_{i}|_{x=X}$,
namely we take the expectation with respect to the pair $Y,\tilde{W}$
with $x$ deterministic, and we replace in the result $x$ by the
random variable $X.$ When this is done, we take the expectation in
the pair $X,Z,$ the only random variables remaining. We can then
rename $Y,\tilde{W}$ into $X,W$ and $X,Z$ into $Y,\tilde{Z}$ in
which $\tilde{Z}$ is a copy of $Z,$ such that $Y,\tilde{Z}$ is
independent of $(X,W).$ We get ( with $y_{j}$ replacing accordingly
$x_{j}$ for the coherence)

\[
(D_{X}\mathbb{E}\sum_{i}V_{i}(Y,\mathcal{L}_{X})W_{i},Z)=\sum_{j}\mathbb{E}_{Y\tilde{Z}}\,\tilde{Z}_{j}\frac{\partial}{\partial y_{j}}\frac{\partial}{\partial m}[\, \mathbb{E}_{XW}\,\sum_{i}V_{i}(X,\mathcal{L}_{X})(Y)W_{i}]=
\]

\[
=\mathbb{E}_{XW}\sum_{i}W_{i}\sum_{j}\mathbb{E}_{Y\tilde{Z}}\tilde{Z}_{j}\frac{\partial}{\partial y_{j}}\frac{\partial}{\partial m}V_{i}(X,\mathcal{L}_{X})(Y)=
\]

\begin{align*}
 & =(W,\sum_{j}\mathbb{E}_{Y\tilde{Z}}\,\tilde{Z}_{j}\frac{\partial}{\partial y_{j}}\frac{\partial}{\partial m}V(X,\mathcal{L}_{X})(Y))=\\
 & =(W,\mathbb{E}_{Y\tilde{Z}}\, D_{y}\frac{\partial}{\partial m}V(X,\mathcal{L}_{X})(Y)\tilde{Z})
\end{align*}

On the other hand 

\[
\mathbb{E}\sum_{j}Z_{j}\frac{\partial V_{i}}{\partial x_{j}}(X,\mathcal{L}_{X})W_{i}=(W,D_{x}V(X,\mathcal{L}_{X})\, Z)
\]
 therefore we conclude that 

\begin{equation}
D\mathcal{U}(X)Z=D_{x}V(X,\mathcal{L}_{X})Z+\mathbb{E}_{Y\tilde{Z}}\, D_{y}\frac{\partial}{\partial m}V(X,\mathcal{L}_{X})(Y)\,\tilde{Z}\label{eq:1-35}
\end{equation}

Suppose we apply this formula to $V(x,m)=D_{x}U(x,m)$ where $U(x,m)$
is scalar, then $\mathcal{U}(X)=D_{x}U(X,\mathcal{L}_{X})$ therefore
(\ref{eq:1-35}) implies 

\begin{equation}
D\mathcal{U}(X)Z=D_{x}^{2}U(X,\mathcal{L}_{X})Z+D_{x}\mathbb{E}_{Y\tilde{Z}}D_{y}\frac{\partial}{\partial m}U(X,\mathcal{L}_{X})(Y).\tilde{Z}\label{eq:1-36}
\end{equation}

In the case $U(x,m)=\dfrac{\partial}{\partial m}F(m)(x),$ then $\mathcal{U}(X)=DF(\mathcal{L}_{X}))$,
hence (\ref{eq:1-36}) implies 

\begin{equation}
D^{2}F(\mathcal{L}_{X})Z=D_{x}^{2}\dfrac{\partial}{\partial m}F(\mathcal{L}_{X})(X)Z+\mathbb{E}_{Y\tilde{Z}}D_{x}D_{y}\frac{\partial^{2}}{\partial m^{2}}F(\mathcal{L}_{X})(X,Y)\tilde{Z}\label{eq:1-37}
\end{equation}

and in the second term $Y,\tilde{Z}$ is an independent copy of $(X,Z)$.

\subsection{APPLICATIONS }

Suppose we apply (\ref{eq:1-37}) to $Z=B(X)\Gamma,$ in wich $\Gamma$
is independent of $X$ and satisfies $\mathbb{E}\Gamma=0,$ then we can take
$\tilde{Z}=B(Y)\tilde{\Gamma}$ , in which $\tilde{\Gamma}$ is a
copy of $\Gamma$with $\tilde{\Gamma}$ independent of $Y$ and the
pair $(Y,\tilde{\Gamma})$ is independent of $X,\Gamma.$ In that
case, the 2nd term vanishes and it remains 

\begin{equation}
D^{2}F(\mathcal{L}_{X})B(X)\Gamma=D_{x}^{2}\dfrac{\partial}{\partial m}F(\mathcal{L}_{X})(X)B(X)\Gamma\label{eq:1-38}
\end{equation}

We apply this formula with $B(x)=\sigma(x)$ and $\Gamma$ a vector
of $n$ independent variables, also independent of $X$ and with mean
$0$ and variance $1.$ We have 

\begin{align*}
(D^{2}F(\mathcal{L}_{X})\sigma(X)\Gamma,\sigma(X)\Gamma) & =\mathbb{E}D_{x}^{2}\dfrac{\partial}{\partial m}F(\mathcal{L}_{X})(X)\sigma(X)\Gamma.\sigma(X)\Gamma=\\
 & =\mathbb{E}\int_{\mathbb{R}^{n}}D_{x}^{2}\dfrac{\partial}{\partial m}F(m)(x)\sigma(x)\Gamma.\sigma(x)\Gamma m(x)dx
\end{align*}

from which it follows easily 

\begin{equation}
(D^{2}F(\mathcal{L}_{X})\sigma(X)\Gamma,\sigma(X)\Gamma)=\int_{\mathbb{R}^{n}}\text{tr (}a(x)D_{x}^{2}\dfrac{\partial}{\partial m}F(m)(x))m(x)dx\label{eq:1-39}
\end{equation}

We next derive another formula. Denote by $e_{k},$ $k=1,\cdots,n$
the vector coordinates in $\mathbb{R}^{n}.$Consider 

\begin{align*}
D^{2}F(\mathcal{L}_{X})e_{k} & =D_{x}^{2}\dfrac{\partial}{\partial m}F(\mathcal{L}_{X})(X)e_{k}+D_{x}\int_{\mathbb{R}^{n}}D_{y}\frac{\partial^{2}}{\partial m^{2}}F(m)(X,y)e_{k}m(y)dy=\\
 & =D_{x}^{2}\dfrac{\partial}{\partial m}F(\mathcal{L}_{X})(X)e_{k}-D_{x}\int_{\mathbb{R}^{n}}\frac{\partial^{2}}{\partial m^{2}}F(m)(X,y)\frac{\partial m(y)}{\partial y_{k}}dy
\end{align*}
 Therefore 

\[
(D^{2}F(\mathcal{L}_{X})e_{k},e_{k})=\int_{\mathbb{R}^{n}}\frac{\partial^{2}}{\partial x_{k}^{2}}\dfrac{\partial}{\partial m}F(m)(x)m(x)dx+\int_{\mathbb{R}^{n}}\frac{\partial^{2}}{\partial m^{2}}F(m)(x,y)\frac{\partial m(x)}{\partial x_{k}}\frac{\partial m(y)}{\partial y_{k}}dy
\]

hence 

\begin{equation}
\sum_{k=1}^{n}(D^{2}F(\mathcal{L}_{X})e_{k},e_{k})=\int_{\mathbb{R}^{n}}\dfrac{\partial}{\partial m}F(m)(x)\Delta m(x)dx+\int_{\mathbb{R}^{n}}\int_{\mathbb{R}^{n}}\frac{\partial^{2}}{\partial m^{2}}F(m)(x,y)Dm(x).Dm(y)dxdy\label{eq:1-40}
\end{equation}

Using (\ref{eq:1-39}) this is also written as 

\begin{equation}
\sum_{k=1}^{n}(D^{2}F(\mathcal{L}_{X})e_{k},e_{k})-(D^{2}F(\mathcal{L}_{X})\Gamma,\Gamma)=\int_{\mathbb{R}^{n}}\int_{\mathbb{R}^{n}}\frac{\partial^{2}}{\partial m^{2}}F(m)(x,y)Dm(x).Dm(y)dxdy\label{eq:1-41}
\end{equation}

We can recover formulas (\ref{eq:1-40}), (\ref{eq:1-41}) by an approach
relying on stochastic differential equations (SDE)

We first consider the solution of the SDE

\begin{align}
dx & =\sigma(x)dw(t)\label{eq:1-6}\\
x(0) & =X\nonumber 
\end{align}

and $\mathcal{L}_{X}=m.$ The Wiener process $w(t)$ is a standard
Wiener process in $\mathbb{R}^{n}$ which is independent of $X.$ Define the
matrix $a(x)=\sigma(x)\sigma^{*}(x).$ The random variable $x(t)$
has a density $m(x,t)$ solution of the Fokker-Planck equation 

\begin{align*}
\frac{\partial m}{\partial t}-\frac{1}{2}\sum_{i,j=1}^{n}\frac{\partial^{2}}{\partial x_{i}\partial x_{j}}(a_{i,j}(x)m) & =0\\
m(x,0)=m(x)
\end{align*}

Therefore 

\begin{align*}
\frac{d}{dt}F(m(t)) & =\int\frac{\partial F}{\partial m}(m(t))(x)\frac{\partial m}{\partial t}(x,t)dx=\\
= & \frac{1}{2}\sum_{i,j=1}^{n}\int a_{i,j}(x)\frac{\partial^{2}}{\partial x_{i}\partial x_{j}}\frac{\partial F}{\partial m}(m(t))(x)m(x,t)dx
\end{align*}
 from which it follows 

\[
\frac{d}{dt}F(m(t))|_{t=0}=\frac{1}{2}\int\text{tr }(a(x)D_{x}^{2}\frac{\partial F}{\partial m}(m)(x))m(x)dx
\]
 but because of (\ref{eq:1-5}) we can state immediately the formula 

\begin{equation}
D_{x}^{2}\frac{\partial F}{\partial m}(m)(x)=D_{x}\partial_{m}F(m)(x)\label{eq:1-7}
\end{equation}

therefore 

\begin{equation}
\frac{d}{dt}F(m(t))|_{t=0}=\frac{1}{2}\int\text{tr }(a(x)D_{x}\partial_{m}F(m)(x))m(x)dx\label{eq:1-8}
\end{equation}

On the other hand, for $t$ small we can approximate the process $x(t)$
by 

\[
x(t)\thicksim X+\sqrt{t}\sigma(X)N
\]
 where $N$ is a standard gaussian independent of $X$ with values
in $\mathbb{R}^{n}.$ We can then write 

\[
F(m(t))\sim F(X+\sqrt{t}\sigma(X)N)
\]

and from formula (\ref{eq:1.2})

\begin{align*}
F(X+\sqrt{t}\sigma(X)N) & =F(X)+(DF(X),\sqrt{t}\sigma(X)N)+\\
+ & \int_{0}^{1}\int_{0}^{1}\lambda D^{2}F(X+\lambda\mu\sqrt{t}\sigma(X)N)(\sqrt{t}\sigma(X)N,\sqrt{t}\sigma(X)N)
\end{align*}
 The 2nd term on the right hand side is $0,$ since $N$ and $X$
are independent. It follows easily 

\[
F(X+\sqrt{t}\sigma(X)N)\sim F(X)+\frac{t}{2}D^{2}F(X)(\sigma(X)N,\sigma(X)N)
\]

which imediately leads to 

\begin{align}
D^{2}F(X)(\sigma(X)N,\sigma(X)N) & =\int\text{tr }(a(x)D_{x}^{2}\frac{\partial F}{\partial m}(m)(x))m(x)dx=\label{eq:1-9}
\end{align}

Consider now the following model 

\[
x(t)=X+\beta b(t)
\]

in which $b(t)$ is a standard Wiener process in $\mathbb{R}^{n},$ and $\beta$
is a constant. We denote by $\mathcal{B}^{t}$ the $\sigma-$algebra
generated by $b(s),s\leq t.$ We consider then the conditional probability
of $x(t)$ given $\mathcal{B}^{t},$ so $\mathcal{L}_{x(t)}^{\mathcal{B}^{t}}$,
which we denote again by $m(t)=m(x,t).$ We have obviously 

\[
m(x,t)=m(x-\beta b(t))
\]
 so the function $m(x,t)$ is random and $\mathcal{B}^{t}$ measurable,
as of course expected from the definition. This time $m(x,t)$ satisfies
a stochastic partial differential equation SPDE, 

\begin{align*}
d_{t}m-\frac{\beta^{2}}{2}\Delta m\, dt+\beta Dm.db & =0\\
m(x,0)=m(x)
\end{align*}

From (\ref{eq:1-1}) we can derive the following Ito's formula 

\begin{align}
dF(m(t)) & =\int\frac{\partial F}{\partial m}(m(t))(\xi)[\frac{\beta^{2}}{2}\Delta m(\xi,t)dt-\beta Dm(\xi,t).db(t)]d\xi\int+\label{eq:1-10}\\
+\frac{\beta^{2}}{2} & \int\int\frac{\partial^{2}F}{\partial m^{2}}(m(t))(\xi,\eta)Dm(\xi,t).Dm(\eta,t)d\xi d\eta dt\nonumber 
\end{align}

from which we get easily 

\begin{align}
\frac{d}{dt}\mathbb{E}F(m(t))|_{t=0} & =\frac{\beta^{2}}{2}\int\frac{\partial F}{\partial m}(m)(\xi)\Delta m(\xi)d\xi+\label{eq:1-11}\\
+ & \frac{\beta^{2}}{2}\int\int\frac{\partial^{2}F}{\partial m^{2}}(m)(\xi,\eta)Dm(\xi).Dm(\eta)d\xi d\eta\nonumber 
\end{align}

On the other hand $F(m(t))=F(X+\beta b(t)).$Using formula (\ref{eq:1-2})
we have 

\begin{align*}
F(X+\beta b(t)) & =F(X)+\beta(DF(X),b(t))+\\
+ & \beta^{2}\,\int_{0}^{1}\int_{0}^{1}\lambda D^{2}F(X+\lambda\mu\beta b(t))(b(t),b(t))
\end{align*}

from which it follows easily that 

\begin{equation}
\frac{d}{dt}\mathbb{E}F(m(t))|_{t=0}=\frac{\beta^{2}}{2}\sum_{k=1}^{n}D^{2}F(X)(e_{k},e_{k})\label{eq:1-12}
\end{equation}

in which $e_{k}$ are the coordinates vectors in $\mathbb{R}^{n}.$ Comparing
with (\ref{eq:1-11}) we have shown the relation 

\begin{equation}
\sum_{k=1}^{n}D^{2}F(X)(e_{k},e_{k})=\int\frac{\partial F}{\partial m}(m)(\xi)\Delta m(\xi)d\xi+\int\int\frac{\partial^{2}F}{\partial m^{2}}(m)(\xi,\eta)Dm(\xi).Dm(\eta)d\xi d\eta\label{eq:1-13}
\end{equation}

Using (\ref{eq:1-9}) we get also 

\begin{equation}
\int\int\frac{\partial^{2}F}{\partial m^{2}}(m)(\xi,\eta)Dm(\xi).Dm(\eta)d\xi d\eta=\sum_{k=1}^{n}D^{2}F(X)(e_{k},e_{k})-D^{2}F(X)(N,N)\label{eq:1-14}
\end{equation}

\begin{rem}
\label{rem1} If $F(\mathcal{L}_{X})=\int_{\mathbb{R}^{n}}F(x)m(x)dx,$ then
it follows from (\ref{eq:1-41}) that $\sum_{k=1}^{n}(D^{2}F(\mathcal{L}_{X})e_{k},e_{k})-(D^{2}F(\mathcal{L}_{X})\Gamma,\Gamma)=0.$ 
\end{rem}

\subsection{APPROACH OF BUCKDAHN et al.}

Remembering the definition of $\partial_{m}F(m)(x),$see (\ref{eq:1-4}),
it is natural to extend it to the 2nd derivative, considering $x$
fixed and thus $\partial_{m}F(m)(x)$ as a function of $m.$ This
is the approach proposed by \cite{BLPR}. Note that $\partial_{m}F(m)(x)$
is a vector in $\mathbb{R}^{n}.$ So consider the component $\partial_{m,i}F(m)(x)$.
We associate to it the function of the random variable $X$ , denoted
$\partial_{m,i}F(X)(x)$. We can define its derivative $D\,\partial_{m,i}F(X)(x)$
as in (\ref{eq:1-11-1}), if of course, the function has a derivative.
But then there exists a new function $\partial_{m}(\partial_{m,i}F(m)(x))(y)$
such that 

\[
\partial_{m}(\partial_{m,i}F(m)(x))(X)=D\,\partial_{m,i}F(X)(x)
\]

We then write 

\begin{equation}
\partial_{m,j}(\partial_{m,i}F(m)(x))(y)=\partial_{m,ij}^{2}F(m)(x,y)\label{eq:1-15}
\end{equation}

This is a new concept. We can from the definitions establish the relation

\begin{equation}
\partial_{m,ij}^{2}F(m)(x,y)=D_{y_{j}}\frac{\partial}{\partial m}(\, D_{x_{i}}\frac{\partial F(m)}{\partial m}(x))(y)\label{eq:1-16}
\end{equation}

We can check that 

\[
\frac{\partial}{\partial m}(\, D_{x_{i}}\frac{\partial F(m)}{\partial m}(x))(y)=D_{x_{i}}\dfrac{\partial^{2}F(m)}{\partial m^{2}}(x,y)
\]
 therefore 

\begin{equation}
\partial_{m,ij}^{2}F(m)(x,y)=D_{y_{j}}D_{x_{i}}\dfrac{\partial^{2}F(m)}{\partial m^{2}}(x,y)\label{eq:1-17}
\end{equation}

So we may write 

\begin{equation}
\partial_{m}^{2}F(m)(x,y)=D_{y}D_{x}\dfrac{\partial^{2}F(m)}{\partial m^{2}}(x,y)=D_{x}D_{y}\dfrac{\partial^{2}F(m)}{\partial m^{2}}(x,y)\label{eq:1-171}
\end{equation}
 and going back to (\ref{eq:1-37}) we obtain 

\[
D^{2}F(\mathcal{L}_{X})Z=D_{x}^{2}\dfrac{\partial}{\partial m}F(\mathcal{L}_{X})(X)Z+\mathbb{E}_{Y\tilde{Z}}(\partial_{m}^{2}F(m)(X,Y)\tilde{Z)}
\]

which gives the relation between the second derivative introduced
by Buckdahn et al. and the second derivative in the Hilbert space
$\mathcal{H}.$ 

In this set up, Taylor's formula is more complex. Consider two probabilities
$m_{0}$ and $m$ to which correspond random variables $X_{0}$ and
$X.$ We use (\ref{eq:1-2}) with the notation of (\cite{BLPR}). 

We have 

\begin{align*}
F(X)-F(X_{0}) & =(DF(X_{0}),X-X_{0})+\int_{0}^{1}\int_{0}^{1}\lambda D^{2}F(X_{0}+\lambda\mu(X-X_{0}))(X-X_{0},X-X_{0})d\lambda d\mu
\end{align*}

Now from formula (\ref{eq:1-37}) we can write 

\begin{equation}
(D^{2}F(\mathcal{L}_{X})Z,Z)=\mathbb{E}D_{x}^{2}\dfrac{\partial}{\partial m}F(\mathcal{L}_{X})(X)Z.Z+\mathbb{E}_{XZ}\, Z.\mathbb{E}_{Y\tilde{Z}}D_{x}D_{y}\frac{\partial^{2}}{\partial m^{2}}F(\mathcal{L}_{X})(X,Y)\tilde{Z}\label{eq:1-172}
\end{equation}

Using this formula and notattion (\cite{BLPR}) we obtain 

\begin{align}
F(m)-F(m_{0}) & =(\,\partial_{m}F(m_{0})(X_{0}),\: X-X_{0})+\label{eq:1-21}
\end{align}
\begin{align*}
+\int_{0}^{1}\int_{0}^{1}\lambda \mathbb{E}_{XX_{0}}(X-X_{0}).\, \mathbb{E}_{YY_{0}}\,\partial_{m}^{2}F(\mathcal{L}_{X_{0}+\lambda\mu(X-X_{0})})(X_{0}+\lambda\mu(X-X_{0}),Y_{0}+\lambda\mu(Y-Y_{0}))(Y-Y_{0})d\lambda d\mu+
\end{align*}
\[
+\int_{0}^{1}\int_{0}^{1}\lambda \mathbb{E}D_{x}\partial_{m}F(\mathcal{L}_{X_{0}+\lambda\mu(X-X_{0})})(X_{0}+\lambda(X-X_{0}))(X-X_{0}).(X-X_{0})d\lambda d\mu
\]

which we write as follows 

\begin{equation}
F(m)-F(m_{0})=\mathbb{E}\,\partial_{m}F(m_{0})(X_{0}).(X-X_{0})+\label{eq:1-22}
\end{equation}

\[
+\frac{1}{2}\mathbb{E}_{XX_{0}}(X-X_{0}).\, \mathbb{E}_{YY_{0}}\,\partial_{m}^{2}F(m_{0})(X_{0},Y_{0})(Y-Y_{0})+\frac{1}{2}\mathbb{E}D_{x}\partial_{m}F(m_{0})(X_{0})(X-X_{0}).(X-X_{0})+R(m_{0};m)
\]
 Buckdahn et al. \cite{BLPR} have shown that $R(m_{0};m)$$\rightarrow0$
with the order of $\mathbb{E}|X-X_{0}|^{3},$ with appropriate assumptions
on the second derivatives. Formulas (\ref{eq:1-21}) and (\ref{eq:1-22})
are substantially different from (\ref{eq:1-1}) or (\ref{eq:1-2}).
This is because of the two variables $m$ and $x$ instead of a single
one $m$ or $X.$

\section{DYNAMIC PROGRAMMING}

\subsection{MEAN-FIELD TYPE CONTROL }

We consider a probability space $\Omega,\mathcal{A},P$ on which are
defined various Wiener processes. We define first a standard Wiener
process $w(t)$ in $\mathbb{R}^{n}$. The classical mean field type control
problem is the following: Let $v(x)$ be a feedback with values in
$\mathbb{R}^{d}$, the corresponding state equation is the Mac Kean-Vlasov
equation 

\begin{align}
dx & =g(x,m_{v(.)},v(x))dt+\sigma(x)dw\label{eq:1.1}\\
x(0) & =x_{0}\nonumber 
\end{align}
in which $m_{v(.)}(x,t)$ is the probability density of the state
$x(t).$ The initial value $x_{0}$ is a random variable independent
of the Wiener process $w(.).$ This density is well defined provided
we assume the invertibility of $a(x)=\sigma(x)\sigma^{*}(x)$ . We
define the 2nd order differential operator 

\[
A\varphi(x)=-\frac{1}{2}\sum_{i,j}a_{ij}(x)\frac{\partial^{2}\varphi(x)}{\partial x_{i}\partial x_{j}}
\]
 and its adjoint 

\[
A^{*}\varphi(x)=-\frac{1}{2}\sum_{i,j}\frac{\partial^{2}(a_{ij}(x)\varphi(x))}{\partial x_{i}\partial x_{j}}
\]
 Next define

\begin{equation}
J(v(.))=\mathbb{E}[\int_{0}^{T}f(x(t),m_{v(.)}(t),v(x(t)))dt+h(x(T),m_{v(.)}(T))]\label{eq:1.2}
\end{equation}
The mean field type control problem is to minimize $J(v(.))$. Note
that the feedback $v(x)$ can also depend on $m.$ It will be indeed
the case for the optimal one. The problem can be easily transformed
into a control problem for a state, which is the probability density
$m_{v(.)}$. It is the solution of the Fokker-Planck equation 

\begin{align}
\frac{\partial m_{v(.)}}{\partial t}+A^{*}m_{v(.)}+\text{div }(g(x,m_{v(.)},v(x))m_{v(.)}(x))\label{eq:2.1}\\
m(x,0)=m_{0}(x)\nonumber 
\end{align}
in which $m_{0}(x)$ is the probability density of the initial value
$x_{0}.$ The objective functional $J(v(.))$ can be written as 

\begin{equation}
J(v(.))=\int_{0}^{T}\int_{\mathbb{R}^{n}}f(x,m_{v(.)}(t),v(x))m_{v(.)}(x,t)dxdt+\int_{\mathbb{R}^{n}}h(x,m_{v(.)}(T))m_{v(.)}(x,T)dx\label{eq:2.2}
\end{equation}

We next use the traditional invariant embedding approach. Define a
family of problems, indexed by initial conditions $m,t$

\begin{align}
\frac{\partial m_{v(.)}}{\partial s}+A^{*}m_{v(.)}+\text{div }(g(x,m_{v(.)},v(x))m_{v(.)}(x))\label{eq:2.3}\\
m(x,t)=m(x)\nonumber 
\end{align}

\begin{equation}
J_{m,t}(v(.))=\int_{t}^{T}\int_{\mathbb{R}^{n}}f(x,m_{v(.)}(s),v(x))m_{v(.)}(x,s)dxds+\int_{\mathbb{R}^{n}}h(x,m_{v(.)}(T))m_{v(.)}(x,T)dx\label{eq:2.4}
\end{equation}

and we set 

\begin{equation}
V(m,t)=\inf_{v(.)}J_{m,t}(v(.))\label{eq:2.5}
\end{equation}

Using the optimality principle one obtains

\begin{align}
\frac{\partial V}{\partial t}-\int_{\mathbb{R}^{n}}\frac{\partial V(m)}{\partial m}(\xi)A^{*}m(\xi)d\xi+\label{eq:2.6}
\end{align}

\[
+\inf_{v}\left(\int_{\mathbb{R}^{n}}f(\xi,m,v(\xi))m(\xi)d\xi\right.\left.-\int_{\mathbb{R}^{n}}\frac{\partial V(m)}{\partial m}(\xi)\text{div }(g(\xi,m,v(\xi))m(\xi))d\xi\right)=0
\]

and 

\begin{equation}
V(m,T)=\int_{\mathbb{R}^{n}}h(x,m)m(x)dx\label{eq:2.7}
\end{equation}

The infimum is attained by minimizing inside the integral. Therefore,
if we define the Hamiltonian $H(x,m,q)$ , with $q\in \mathbb{R}^{n}$ by 

\begin{equation}
H(x,m,q)=\inf_{v}(f(x,m,v)+q.g(x,m,v))\label{eq:2.8-2}
\end{equation}

we can write (\ref{eq:2.6}) as follows 

\begin{equation}
\frac{\partial V}{\partial t}-\int_{\mathbb{R}^{n}}A\,\frac{\partial V(m)}{\partial m}(\xi)m(\xi)d\xi+\int_{\mathbb{R}^{n}}H(\xi,m,D\frac{\partial V(m)}{\partial m}(\xi))m(\xi)d\xi=0\label{eq:2-9}
\end{equation}

So (\ref{eq:2-9}), (\ref{eq:2.7}) is the Bellman equation for the
mean-field type control problem (\ref{eq:2.1}),(\ref{eq:2.2}), with
derivatives with respect to densities. 

The next step is to write the corresponding equation for $V(X,t)$
with the argument $X$ in $\mathcal{H}$ . We use (\ref{eq:1-3})
and (\ref{eq:1-39}). 

Then (\ref{eq:2-9}) becomes 

\begin{align}
\frac{\partial V}{\partial t}+\frac{1}{2}D^{2}V(X)(\sigma(X)\Gamma,\sigma(X)\Gamma)+\mathbb{E}H(X,\mathcal{L}_{X},DV(X)) & =0\label{eq:2-12}\\
V(X,T)=\mathbb{E}\, h(X,\mathcal{L}_{X})\nonumber 
\end{align}

\begin{example}
\label{exa2} If we take $\sigma=0,\, g(x,m,v)=v,\, f(x,m,v)=\frac{1}{2}|v|^{2},$we
get the Eikonal equation 
\[
\frac{\partial V}{\partial t}-\frac{1}{2}||DV(X)||^{2}=0
\]

\end{example}

\subsection{STOCHASTIC MEAN FIELD TYPE CONTROL}

\subsubsection{\label{sub:PRELIMINARIES}PRELIMINARIES}

If we look at the formulation (\ref{eq:2.1}), (\ref{eq:2.2}) of
the mean field type control problem, it is a deterministic problem,
although at the origin it was a stochastic one, see (\ref{eq:1.1}),
(\ref{eq:1.2}). We now consider a stochastic version of (\ref{eq:2.1}),
(\ref{eq:2.2}) or a doubly stochastic version of (\ref{eq:1.1}),
(\ref{eq:1.2}). Let us begin with this one. Assume there is a 2nd
standard Wiener process $b(t)$ with values in $\mathbb{R}^{n}$; $b(t)$ and
$w(t)$ are independent and independent of $x_{0}$. We set $\mathcal{B}^{t}$=$\sigma(b(s),s\leq t)$
and $\mathcal{F}$$^{t}$=$\sigma(x_{0},b(s),w(s),\: s\leq t)$ .
The control $v(x,t)$ at time $t$ is a feedback, but not deterministic.
It is a random variable, $\mathcal{B}^{t}$- measurable. We consider
the stochastic Mc Kean-Vlasov equation 

\begin{align}
dx & =g(x,m_{v(.)}(t),v(x,t))dt+\sigma(x)dw+\beta db(t)\label{eq:3.1}\\
x(0) & =x_{0}\nonumber 
\end{align}
in which $m_{v(.)}(t)$ represents the conditional probability density
of $x(t)$, given the $\sigma-$algebra $\mathcal{B}^{t}$ . We want
to minimize the objective functional 

\begin{equation}
J(v(.))=\mathbb{E}[\int_{0}^{T}f(x(t),m_{v(.)}(t),v(x(t),t))dt+h(x(T),m_{v(.)}(T))]\label{eq:3.2}
\end{equation}

\subsubsection{CONDITIONAL PROBABILITY}

If we define 

\[
y(t)=x(t)-\beta b(t)
\]
 then the process $y(t)$ satisfies the equation 

\begin{align}
dy & =g(y(t)+\beta b(t),m_{v(.)}(t),v(y(t)+\beta b(t),t))dt+\sigma(y(t)+\beta b(t))dw\label{eq:3.1-1}\\
y(0) & =x_{0}\nonumber 
\end{align}
If we fix $b(s),s\leq t$ then the conditional probability of $y(t)$
is simply the probability density arising from the Wiener process
$w(t),$in view of the independence of $w(t)$ and $b(t).$ It is
the function $p(y,t)$ solution of 

\begin{align*}
\frac{\partial p}{\partial t}-\frac{1}{2}\sum_{ij}(a_{ij}(y+\beta b(t))\frac{\partial^{2}p}{\partial y_{i}\partial y_{j}})+
\end{align*}

\[
+\text{div }\left(g(y(t)+\beta b(t),m_{v(.)}(t),v(y(t)+\beta b(t),t))p\right)=0
\]
\[
p(y,0)=m_{0}(y)
\]

The conditional probability density of $x(t)$ given $\mathcal{B}^{t}$
is the function 

\[
m(x,t)=p(x-\beta b(t),t)
\]

hence 
\[
\partial_{t}m=(\frac{\partial p}{\partial t}+\frac{1}{2}\beta^{2}\Delta p)(x-\beta b(t),t)dt-\beta Dp(x-\beta b(t),t)db(t)
\]

which yields 

\begin{equation}
\partial_{t}m+(A^{*}m-\frac{1}{2}\beta^{2}\Delta m+\text{div}(g(x,m(t),v(x,t))m))dt+\beta Dm.db(t)=0\label{eq:3.2-1}
\end{equation}
\[
m(x,0)=m_{0}(x)
\]

and the objective functional (\ref{eq:3.2}) can be written as 

\begin{equation}
J(v(.))=\mathbb{E}[\int_{0}^{T}\int_{\mathbb{R}^{n}}f(x,m(t),v(x,t))m(x,t)dxdt+\int_{\mathbb{R}^{n}}h(x,m(T))m(x,T)dx]\label{eq:3.3}
\end{equation}
The problem becomes a stochastic control problem for a distributed
parameter system. Using invariant embedding, we consider the family
of problems indexed by $m,t$

\begin{equation}
\partial_{s}m+(A^{*}m-\frac{1}{2}\beta^{2}\Delta m+\text{div}(g(x,m(s),v(x,s))m))ds+\beta Dm.db(s)=0\label{eq:3.4}
\end{equation}

\[
m(x,t)=m(x)
\]

and 

\begin{equation}
J_{m,t}(v(.))=\mathbb{E}[\int_{t}^{T}\int_{\mathbb{R}^{n}}f(x,m(s),v(x,s))m(x,s)dxdt+\int_{\mathbb{R}^{n}}h(x,m(T))m(x,T)dx]\label{eq:3.5}
\end{equation}

Set 

\begin{equation}
V(m,t)=\inf_{v(.)}J_{m,t}(v(.))\label{eq:3.51}
\end{equation}
 then $V(m,t)$ satisfies the Dynamic Programming equation 

\begin{align}
\frac{\partial V}{\partial t}-\int_{\mathbb{R}^{n}}\frac{\partial V(m,t)}{\partial m}(\xi)(A^{*}m(\xi)-\frac{1}{2}\beta^{2}\Delta m(\xi))d\xi & +\frac{1}{2}\beta^{2}\int_{\mathbb{R}^{n}}\int_{\mathbb{R}^{n}}\frac{\partial^{2}V(m,t)}{\partial m^{2}}(\xi,\eta)Dm(\xi)Dm(\eta)d\xi d\eta+\label{eq:2.6-1}
\end{align}

\[
+\inf_{v}\left(\int_{\mathbb{R}^{n}}f(\xi,m,v(\xi))m(\xi)d\xi\right.\left.-\int_{\mathbb{R}^{n}}\frac{\partial V(m,t)}{\partial m}(\xi)\text{div }(g(\xi,m,v(\xi))m(\xi))d\xi\right)=0
\]
\[
V(m,T)=\int_{\mathbb{R}^{n}}h(x,m)m(x)dx
\]

which we rewrite as follows 

\begin{equation}
\frac{\partial V}{\partial t}-\int_{\mathbb{R}^{n}}(A-\frac{1}{2}\beta^{2}\Delta)\frac{\partial V(m,t)}{\partial m}(\xi)m(\xi)d\xi+\frac{1}{2}\beta^{2}\int_{\mathbb{R}^{n}}\int_{\mathbb{R}^{n}}\frac{\partial^{2}V(m,t)}{\partial m^{2}}(\xi,\eta)Dm(\xi)Dm(\eta)d\xi d\eta+\label{eq:2.7-2}
\end{equation}
\[
+\int_{\mathbb{R}^{n}}H(\xi,m,D\frac{\partial V(m)}{\partial m}(\xi))m(\xi)d\xi=0
\]
 Let us write the equation for the corresponding function $V(X,t).$
We use now (\ref{eq:1-41}) and we obtain, after cancellation of the
term $\frac{1}{2}\beta^{2}D^{2}V(X)(\Gamma,\Gamma)$ 

\begin{equation}
\frac{\partial V}{\partial t}+\frac{1}{2}D^{2}V(X)(\sigma(X)\Gamma,\sigma(X)\Gamma)+\frac{1}{2}\beta^{2}\sum_{k=1}^{n}D^{2}V(X)(e_{k},e_{k})+\mathbb{E}H(X,\mathcal{L}_{X},DV(X))=0\label{eq:2-8}
\end{equation}
\[
V(X,T)=\mathbb{E}h(X,\mathcal{L}_{X})
\]
 This write up is more condensed than (\ref{eq:2.6-1}).

\section{FIRST TYPE OF MASTER EQUATION}

\subsection{THE PROBLEM OF BUCKDAHN et al.}

In their paper \cite{BLPR}, the authors consider the following problem.
It is not a control problem. As in (\ref{eq:2.3}) , let $m(x,s),s\geq t$
solution of 

\begin{align}
\frac{\partial m}{\partial s}+A^{*}m+\text{div }(g(x,m(s))m(x,s))\label{eq:2.3-1}\\
m(x,t)=m(x)\nonumber 
\end{align}

and we consider the associated diffusion 

\begin{align}
dx & =g(x(s),m(s))ds+\sigma(x(s))dw(s)\label{eq:4.1-1}\\
 & x(t)=x\nonumber 
\end{align}

Of course, if the initial condition $x(t)$ were replaced with a random
variable independant of the Wiener process, with probability density
$m(x),$ then $m(x,s)$ would be simply the probability density of
$x(s).$Since it is not the case, the pair ( \ref{eq:2.3-1}), (\ref{eq:4.1-1})
is a system. The authors are interested in the functional 

\begin{equation}
U(x,m,t)=\mathbb{E}h(x(T),m(T))\label{eq:4.2-1}
\end{equation}

Note indeed that the solution of (\ref{eq:2.3-1}) can be written
as $m_{m,t}(x,s)$ and the solution of (\ref{eq:4.1-1}) as $x_{m,x,t}(s)$,
to emphasize the dependence in the initial conditions, and this justifies
the notation (\ref{eq:4.2-1}). We will write the formulation of the
equation satisfied by $U$ with the 3 approaches explained above.

\subsection{VARIOUS FORMULATIONS }

We begin by considering that $m$ is evolving in $L^{2}(\mathbb{R}^{n}).$
We can see formally , using the Markov property of the pair $m(s),x(s)$
that 

\begin{equation}
\frac{\partial U}{\partial t}-\int_{\mathbb{R}^{n}}\frac{\partial U(x,m,t)}{\partial m}(\xi)A^{*}m(\xi)d\xi-\int_{\mathbb{R}^{n}}\frac{\partial U(x,m,t)}{\partial m}(\xi)\text{div }(g(\xi,m)m(\xi))d\xi+\label{eq:4-3}
\end{equation}

\[
+D_{x}U.g(x,m)-AU(x,m,t)=0
\]
\[
U(x,m,T)=h(x,m)
\]
 We write it as 

\begin{equation}
\frac{\partial U}{\partial t}-\int_{\mathbb{R}^{n}}A_{\xi}\frac{\partial U(x,m,t)}{\partial m}(\xi)m(\xi)d\xi+\int_{\mathbb{R}^{n}}D_{\xi}\frac{\partial U(x,m,t)}{\partial m}(\xi).g(\xi,m)m(\xi)d\xi+\label{eq:4-4}
\end{equation}
\[
+D_{x}U.g(x,m)-A_{x}U(x,m,t)=0
\]

Because of the two variables $x,\xi$ it is important to clarify on
which variables the differential operators act. For that purpose we
have used the notation $D_{x},D_{\xi}$ and $A_{x},$$A_{\xi}.$ Using
formulas (\ref{eq:1-5}) and (\ref{eq:1-7}) we have 

\begin{equation}
\frac{\partial U}{\partial t}+\frac{1}{2}\int_{\mathbb{R}^{n}}\text{tr}(a(\xi)D_{\xi}\partial_{m}U(x,m,t)(\xi))m(\xi)d\xi+\int_{\mathbb{R}^{n}}\partial_{m}U(x,m,t)(\xi).g(\xi,m)m(\xi)d\xi+\label{eq:4-5}
\end{equation}
\[
+D_{x}U.g(x,m)-A_{x}U(x,m,t)=0
\]
\[
U(x,m,T)=h(x,m)
\]
 This is the equation studied rigorously by Buckdahn et al \cite{BLPR}.
The authors notice that there is no second derivative in $m,$although
there are second derivatives in $x.$ This is explained by the fact
that the equation in $m$ is deterministic, whereas the equation in
$x$ is random. Let us now write the equation in the space $\mathcal{H}$.
We consider $U(x,X,t)=U(x,\mathcal{L}_{X},t).$

From (\ref{eq:1-3}) and (\ref{eq:1-9}) we obtain 

\begin{equation}
\frac{\partial U}{\partial t}+\frac{1}{2}D_{X}^{2}U(x,X,t)(\sigma(X)\Gamma,\sigma(X)\Gamma)+(D_{X}U(x,X,t),g(X,\mathcal{L}_{X},t))+\label{eq:4-5-1}
\end{equation}
\[
+D_{x}U(x,X,t).g(x,\mathcal{L}_{X},t)-A_{x}U(x,X,t)=0
\]
\[
U(x,X,T)=h(x,\mathcal{L}_{X})
\]

In this set up, we obtain a more symmetric formulation. There is no
mixed derivatives and there is a second order derivative in $X.$
\begin{rem}
\label{rem3} From the methods described above, it is possible to
extend the equation (\ref{eq:4-5-1}) to nonlinear situations. 
\end{rem}

\section{OBTAINING THE MASTER EQUATION BY DIFFERENTIATING BELLMAN EQUATION }

\subsection{THE MEAN FIELD TYPE CONTROL CASE.}

We first go back to (\ref{eq:2.8-2}) and assume that the minimum
is attained in a unique point which we denote by $\hat{v}(x,m,q).$
We then set 

\begin{equation}
G(x,m,q)=g(x,m,\hat{v}(x,m,q))\label{eq:5-1}
\end{equation}

We have also

\begin{equation}
G(x,m,q)=D_{q}H(x,m,q)\label{eq:5-2}
\end{equation}

We then consider the value function $V(m,t)$ defined in (\ref{eq:3.51})
and set 

\begin{equation}
U(x,m,t)=\frac{\partial V(m,t)}{\partial m}(x)\label{eq:2.8}
\end{equation}

we write (\ref{eq:2-9}) as 

\begin{equation}
\frac{\partial V}{\partial t}-\int_{\mathbb{R}^{n}}AU(\xi,m,t)m(\xi)d\xi+\int_{\mathbb{R}^{n}}H(\xi,m,DU(\xi,m,t))m(\xi)d\xi=0\label{eq:2.9}
\end{equation}
 We next differentiate (\ref{eq:2.9}) in $m$ . We note that 

\[
\frac{\partial}{\partial m}[\int_{\mathbb{R}^{n}}H(\xi,m,DU(\xi,m,t))m(\xi)d\xi](x)=H(x,m,DU(x))+
\]

\[
+\int_{\mathbb{R}^{n}}\frac{\partial}{\partial m}H(\xi,m,DU(\xi))(x)m(\xi)d\xi+\int_{\mathbb{R}^{n}}G(\xi,m,DU(\xi))m(\xi)D_{\xi}\frac{\partial}{\partial m}U(\xi,m,t)(x)d\xi
\]

Hence, using 

\[
\frac{\partial}{\partial m}U(\xi,m,t)(x)=\frac{\partial}{\partial m}U(x,m,t)(\xi)=\frac{\partial^{2}V(m,t)}{\partial m^{2}}(x,\xi)
\]
 we obtain 

\begin{align}
-\frac{\partial U}{\partial t}+AU+\int_{\mathbb{R}^{n}}\frac{\partial}{\partial m}U(x,m,t)(\xi)(A^{*}m(\xi)+\text{div (}G(\xi,m,DU(\xi))m(\xi))d\xi & =\label{eq:2.10}\\
H(x,m,DU(x))+\int_{\mathbb{R}^{n}}\frac{\partial}{\partial m}H(\xi,m,DU(\xi))(x)m(\xi)d\xi\nonumber 
\end{align}
\begin{equation}
U(x,m,T)=h(x,m)+\int_{\mathbb{R}^{n}}\frac{\partial}{\partial m}h(\xi,m)(x)m(\xi)d\xi\label{eq:2.11}
\end{equation}

We write this equation as 

\begin{align}
-\frac{\partial U}{\partial t}+A_{x}U+\int_{\mathbb{R}^{n}}(A_{\xi}\frac{\partial}{\partial m}U(x,m,t)(\xi)-D_{\xi}\frac{\partial}{\partial m}U(x,m,t)(\xi).G(\xi,m,DU(\xi)))m(\xi)d\xi & =\label{eq:2.10-2}\\
H(x,m,DU(x))+\int_{\mathbb{R}^{n}}\frac{\partial}{\partial m}H(\xi,m,DU(\xi))(x)m(\xi)d\xi\nonumber 
\end{align}

\[
U(x,m,T)=h(x,m)+\int_{\mathbb{R}^{n}}\frac{\partial}{\partial m}h(\xi,m)(x)m(\xi)d\xi
\]

\subsection{MASTER EQUATION IN $\mathcal{H}$}

We then introduce $\mathcal{U}(X,t)=DV(\mathcal{L}_{X},t).$ So $\mathcal{U}$
maps $\mathcal{H}\times[0,T]$ into $\mathcal{H}.$ To simplify a
little we take $\sigma(x)=\sigma I$, then Bellman equation (\ref{eq:2-12})
reads 

\begin{align}
\frac{\partial V}{\partial t}+\frac{\sigma^{2}}{2}D^{2}V(X)(\Gamma,\Gamma)+\mathbb{E}\, H(X,\mathcal{L}_{X},\mathcal{U}(X)) & =0\label{eq:5-4}\\
V(X,T)=\mathbb{E}\, h(X,\mathcal{L}_{X})\nonumber 
\end{align}

We proceed with the formal differentiation of (\ref{eq:5-4}) . First
$D\mathcal{U}(X):\mathcal{H}\times[0,T]\rightarrow\mathcal{L}(\mathcal{H};\mathcal{H}).$
So $D\mathcal{U}(X)$ is a random $n\times n$ matrix . We obtain 

\begin{equation}
D\, \mathbb{E}\, H(X,\mathcal{L}_{X},\mathcal{U}(X))=D_{x}H(X,\mathcal{L}_{X},\mathcal{U}(X))+(D\mathcal{U}(X))^{*}G(X,\mathcal{L}_{X},\mathcal{U}(X))+\label{eq:5-6-1}
\end{equation}

\[
+\mathbb{E}_{Y,\mathcal{U}(Y)}D_{x}\frac{\partial H}{\partial m}(Y,\mathcal{L}_{X},\mathcal{U}(Y))(X)
\]

in which $Y,\mathcal{U}(Y)$ is an independent copy of $X,$$\mathcal{U}(X).$
But $D\mathcal{U}(X)=D^{2}V(\mathcal{L}_{X},t)$ is self-adjoint ,
so we obtain that $\mathcal{U}(X,t)$ satisfies the equation 

\begin{equation}
\frac{\partial}{\partial t}\mathcal{U}+\frac{\sigma^{2}}{2}D^{2}\mathcal{U}(\Gamma,\Gamma)+D\mathcal{U}(X)G(X,\mathcal{L}_{X},\mathcal{U}(X))+\label{eq:5-61}
\end{equation}

\[
D_{x}H(X,\mathcal{L}_{X},\mathcal{U}(X))+\mathbb{E}_{Y,\mathcal{U}(Y)}D_{x}\frac{\partial H}{\partial m}(Y,\mathcal{L}_{X},\mathcal{U}(Y))(X)=0
\]

\[
\mathcal{U}(X,T)=D_{x}h(X,\mathcal{L}_{X})+\mathbb{E}_{Y}D_{x}\frac{\partial h}{\partial m}(Y,\mathcal{L}_{X})(X)
\]

We also can obtain this equation from (\ref{eq:2.10-2}). For $\sigma(x)=\sigma I,$
this equation writes 

\begin{equation}
\frac{\partial U}{\partial t}+\frac{\sigma^{2}}{2}\Delta_{x}U+\frac{\sigma^{2}}{2}\int_{\mathbb{R}^{n}}\Delta_{\xi}\frac{\partial}{\partial m}U(x,m,t)(\xi)m(\xi)d\xi+\int_{\mathbb{R}^{n}}D_{\xi}\frac{\partial}{\partial m}U(x,m,t)(\xi).G(\xi,m,DU(\xi)))m(\xi)d\xi+\label{eq:5-6}
\end{equation}

\[
+H(x,m,D_{x}U(x))+\int_{\mathbb{R}^{n}}\frac{\partial}{\partial m}H(\xi,m,D_{\xi}U(\xi))(x)m(\xi)d\xi=0
\]
 We differentiate in $x,$so we get 

\begin{equation}
\frac{\partial}{\partial t}D_{x}U+\frac{\sigma^{2}}{2}\Delta_{x}D_{x}U+\frac{\sigma^{2}}{2}D_{x}\int_{\mathbb{R}^{n}}\Delta_{\xi}\frac{\partial}{\partial m}U(x,m,t)(\xi)m(\xi)d\xi+D_{x}\int_{\mathbb{R}^{n}}D_{\xi}\frac{\partial}{\partial m}U(x,m,t)(\xi).G(\xi,m,DU(\xi)))m(\xi)d\xi+\label{eq:5-7}
\end{equation}

\[
+D_{x}H(x,m,D_{x}U(x))+D_{x}^{2}U(x,m,t)G(x,m,D_{x}U)+D_{x}\int_{\mathbb{R}^{n}}\frac{\partial}{\partial m}H(\xi,m,D_{\xi}U(\xi))(x)m(\xi)d\xi=0
\]
 We recall that $D_{x}U(X,m,t)=DV(X,t)=\mathcal{U}(X,t).$ Now from
formula (\ref{eq:1-39}) we have, skipping the argument $t$ 

\[
D^{2}V(X)(\Gamma,\Gamma)=\int_{\mathbb{R}^{n}}\Delta_{x}U(x,m)m(x)dx
\]
 If we differentiate in $X$ the left hand side, we get $D^{2}\mathcal{U}(X)(\Gamma,\Gamma).$
But the right hand side is a functional of $m.$ So we can obtain
the same result , by first differentiating the right hand side in
$m,$then the result in $x,$and finally replacing $x$ by $X.$ Calling
$F(m)$ the right hand side, we first have ( using the symmetry property)

\begin{align*}
\frac{\partial F(m)}{\partial m}(x) & =\Delta_{x}U(x,m)+\int_{\mathbb{R}^{n}}\Delta_{\xi}\frac{\partial}{\partial m}U(\xi,m)(x)m(\xi)d\xi\\
 & =\Delta_{x}U(x,m)+\int_{\mathbb{R}^{n}}\Delta_{\xi}\frac{\partial}{\partial m}U(x,m)(\xi)m(\xi)d\xi
\end{align*}
 Hence
\[
D_{x}\frac{\partial F(m)}{\partial m}(x)=\Delta_{x}D_{x}U(x,m)+D_{x}\int_{\mathbb{R}^{n}}\Delta_{\xi}\frac{\partial}{\partial m}U(x,m)(\xi)m(\xi)d\xi
\]

Replacing $x$ by $X,$ we conclude 

\begin{equation}
D^{2}\mathcal{U}(X)(\Gamma,\Gamma)=\Delta_{x}D_{x}U(X,\mathcal{L}_{X})+D_{x}\int_{\mathbb{R}^{n}}\Delta_{\xi}\frac{\partial}{\partial m}U(X,\mathcal{L}_{X})(\xi)m(\xi)d\xi\label{eq:5-8}
\end{equation}

We next have 

\begin{equation}
D_{x}H(X,m,D_{x}U(X))=D_{x}H(X,\mathcal{L}_{X},\mathcal{U}(X))\label{eq:5-9}
\end{equation}

and 

\begin{equation}
D_{x}\int_{\mathbb{R}^{n}}\frac{\partial}{\partial m}H(\xi,m,D_{\xi}U(\xi))(X)m(\xi)d\xi=\mathbb{E}_{Y,\mathcal{U}(Y)}D_{x}\frac{\partial H}{\partial m}(Y,\mathcal{L}_{X},\mathcal{U}(Y))(X)\label{eq:5-10}
\end{equation}

Finally it remains to check that 

\begin{align}
D_{x}\int_{\mathbb{R}^{n}}D_{\xi}\frac{\partial}{\partial m}U(X,m,t)(\xi).G(\xi,m,D_{\xi}U(\xi)))m(\xi)d\xi+D_{x}^{2}U(X,m,t)G(X,m,D_{x}U(X)) & =\label{eq:5-11}\\
D\mathcal{U}(X,t)G(X,\mathcal{L}_{X},\mathcal{U}(X))\nonumber 
\end{align}

Using formula (\ref{eq:1-37}) with $Z=G(X,\mathcal{L}_{X},\mathcal{U}(X))$
we obtain

\[
D\mathcal{U}(X,t)G(X,\mathcal{L}_{X},\mathcal{U}(X))=D_{x}^{2}U(X,\mathcal{L}_{X},t)G(X,\mathcal{L}_{X},\mathcal{U}(X))+
\]

\[
+\int_{\mathbb{R}^{n}}D_{x}D_{\xi}\frac{\partial}{\partial m}U(X,\mathcal{L}_{X},t)(\xi)G(\xi,\mathcal{L}_{X},D_{\xi}U(\xi))m(\xi)d\xi
\]

which is exactly (\ref{eq:5-11}).

\subsection{\label{sub:INTERPRETATION-OF-THE}INTERPRETATION OF THE MASTER EQUATION}

We can interpret $U(x,m,t)$ as uncoupling the system of HJB-FP equations
of the mean-field type control problem, see \cite{BFY1}. We recall
briefly the idea.

The probability density, corresponding to the optimal feedback control
is given by 

\begin{align*}
\frac{\partial m}{\partial t}+A^{*}m+\text{div }(G(x,m,DU)m(x)) & =0\\
m(x,0)=m_{0}(x)
\end{align*}

Define $u(x,t)=U(x,m(t),t),$ then clearly, from (\ref{eq:2.10}),(\ref{eq:2.11})
we obtain 

\begin{align}
-\frac{\partial u}{\partial t}+Au & =H(x,m,Du(x))+\int_{\mathbb{R}^{n}}\frac{\partial}{\partial m}H(\xi,m,Du(\xi))(x)m(\xi)d\xi\label{eq:2.12}\\
u(x,T) & =h(x,m)+\int_{\mathbb{R}^{n}}\frac{\partial}{\partial m}h(\xi,m)(x)m(\xi)d\xi\label{eq:2.13}
\end{align}
 which together with the F-P equation

\begin{align}
\frac{\partial m}{\partial t}+A^{*}m+\text{div }(G(x,m,Du)m(x)) & =0\label{eq:2.14}\\
m(x,0)=m_{0}(x)\label{eq:2.15}
\end{align}
 form the system of coupled HJB- FP equations of the classical mean
field type control problem, see \cite{BFY}.

\subsection{STOCHASTIC MEAN FIELD TYPE CONTROL }

We go back to (\ref{eq:2.7-2}). We define again $U(x,m,t)=\frac{\partial V(m,t)}{\partial m}(x)$,
so we can write 

\begin{equation}
U(x,m,T)=h(x,m)+\int_{\mathbb{R}^{n}}\frac{\partial h(\xi,m)}{\partial m}(x)m(\xi)d\xi\label{eq:2.6-2}
\end{equation}

and from (\ref{eq:2.7-2}) we obtain

\begin{align}
\frac{\partial V}{\partial t}-\int_{\mathbb{R}^{n}}(AU-\frac{1}{2}\beta^{2}\Delta U)(\xi,m,t)m(\xi,t)d\xi & +\label{eq:2.7-1}
\end{align}

\[
\frac{1}{2}\beta^{2}\int_{\mathbb{R}^{n}}\int_{\mathbb{R}^{n}}\frac{\partial U(\xi,m,t)}{\partial m}(\eta)Dm(\xi)Dm(\eta)d\xi d\eta+\int_{\mathbb{R}^{n}}H(\xi,m,D_{\xi}U)m(\xi)dx=0
\]

We then differentiate this equation in $m$, to get an equation for
$U.$ We note that 

\begin{align*}
\frac{\partial}{\partial m}\left(\int_{\mathbb{R}^{n}}\int_{\mathbb{R}^{n}}\frac{\partial U(\xi,m,t)}{\partial m}(\eta)Dm(\xi)Dm(\eta)d\xi d\eta\right)(x) & =\int_{\mathbb{R}^{n}}\int_{\mathbb{R}^{n}}\frac{\partial^{2}U(x,m,t)}{\partial m^{2}}(\xi,\eta)Dm(\xi)Dm(\eta)d\xi d\eta-\\
- & 2\text{div}\left(\int_{\mathbb{R}^{n}}\frac{\partial U(x,m,t)}{\partial m}(\eta)Dm(\eta)d\eta\right)
\end{align*}

and 

\begin{align*}
\frac{\partial}{\partial m}\left(\int_{\mathbb{R}^{n}}H(\xi,m,DU)m(\xi)dx\right) & =H(x,m,DU(x))+\int_{\mathbb{R}^{n}}\frac{\partial}{\partial m}H(\xi,m,D_{\xi}U(\xi))(x)m(\xi)d\xi+\\
+ & \int_{\mathbb{R}^{n}}\frac{\partial}{\partial m}(D_{\xi}U(\xi,m,t))(x)G(\xi,m,DU(\xi))m(\xi)d\xi
\end{align*}
 Next 

\[
\int_{\mathbb{R}^{n}}\frac{\partial}{\partial m}(D_{\xi}U(\xi,m,t))(x)G(\xi,m,DU(\xi))m(\xi)d\xi=\int_{\mathbb{R}^{n}}D_{\xi}(\frac{\partial}{\partial m}U(\xi,m,t)(x))G(\xi,m,DU(\xi))m(\xi)d\xi=
\]
\[
-\int_{\mathbb{R}^{n}}\frac{\partial}{\partial m}U(\xi,m,t)(x)\:\text{div}(G(\xi,m,DU(\xi))m(\xi))d\xi=-\int_{\mathbb{R}^{n}}\frac{\partial}{\partial m}U(x,m,t)(\xi)\:\text{div}(G(\xi,m,DU(\xi))m(\xi))d\xi
\]
 Collecting results, we obtain the Master equation 

\begin{align}
-\frac{\partial U}{\partial t}+AU-\frac{1}{2}\beta^{2}\Delta U+\label{eq:2.8-1}
\end{align}

\[
+\int_{\mathbb{R}^{n}}(A_{\xi}\frac{\partial}{\partial m}U(x,m,t)(\xi)-D_{\xi}\frac{\partial}{\partial m}U(x,m,t)(\xi).G(\xi,m,DU(\xi)))m(\xi)d\xi-\frac{1}{2}\beta^{2}\int_{\mathbb{R}^{n}}\frac{\partial}{\partial m}U(x,m,t)(\xi)\Delta_{\xi}m(\xi)d\xi
\]

\[
-\frac{1}{2}\beta^{2}\int_{\mathbb{R}^{n}}\int_{\mathbb{R}^{n}}\frac{\partial^{2}U(x,m,t)}{\partial m^{2}}(\xi,\eta)Dm(\xi)Dm(\eta)d\xi d\eta+\beta^{2}\text{div}\left(\int_{\mathbb{R}^{n}}\frac{\partial U(x,m,t)}{\partial m}(\xi)Dm(\xi)d\xi\right)
\]

\[
=H(x,m,DU(x))+\int_{\mathbb{R}^{n}}\frac{\partial}{\partial m}H(\xi,m,DU(\xi))(x)m(\xi)d\xi
\]

\begin{equation}
U(x,m,T)=h(x,m)+\int_{\mathbb{R}^{n}}\frac{\partial h(\xi,m)}{\partial m}(x)m(\xi)d\xi\label{eq:2.9-1}
\end{equation}

We next write the Master equation in the space $\mathcal{H}$ , for
$\mathcal{U}(X)=DV(X)=D_{x}U(X,\mathcal{L}_{X})$. Recall that we
take $\sigma(x)=I.$ From (\ref{eq:2-8}) and (\ref{eq:5-61}) we
get 

\begin{equation}
\frac{\partial}{\partial t}\mathcal{U}+\frac{\sigma^{2}}{2}D^{2}\mathcal{U}(\Gamma,\Gamma)+\frac{1}{2}\beta^{2}\sum_{k=1}^{n}D^{2}\mathcal{U}(X)(e_{k},e_{k})+D\mathcal{U}(X)G(X,\mathcal{L}_{X},\mathcal{U}(X))+\label{eq:5-61-1}
\end{equation}

\[
D_{x}H(X,\mathcal{L}_{X},\mathcal{U}(X))+\mathbb{E}_{Y,\mathcal{U}(Y)}D_{x}\frac{\partial H}{\partial m}(Y,\mathcal{L}_{X},\mathcal{U}(Y))(X)=0
\]

\[
\mathcal{U}(X,T)=D_{x}h(X,\mathcal{L}_{X})+\mathbb{E}_{Y}D_{x}\frac{\partial h}{\partial m}(Y,\mathcal{L}_{X})(X)
\]

If we use the Master equation (\ref{eq:2.8-1}), (\ref{eq:2.9-1})
, then from the calculations done in the case $\beta=0,$ what remains
to be proven is 

\begin{equation}
\sum_{k=1}^{n}D^{2}\mathcal{U}(X)(e_{k},e_{k})=D_{x}\Delta_{x}U(X,\mathcal{L}_{X})+D_{x}\int_{\mathbb{R}^{n}}\frac{\partial}{\partial m}U(X,m,t)(\xi)\Delta_{\xi}m(\xi)d\xi+\label{eq:5-62}
\end{equation}

\[
+D_{x}\int_{\mathbb{R}^{n}}\int_{\mathbb{R}^{n}}\frac{\partial^{2}U(X,m,t)}{\partial m^{2}}(\xi,\eta)Dm(\xi)Dm(\eta)d\xi d\eta-2D_{x}\text{div}\left(\int_{\mathbb{R}^{n}}\frac{\partial U(X,m,t)}{\partial m}(\xi)Dm(\xi)d\xi\right)
\]

From formula (\ref{eq:1-40}) we can write 

\[
\sum_{k=1}^{n}(D^{2}V(\mathcal{L}_{X})e_{k},e_{k})=\int_{\mathbb{R}^{n}}\dfrac{\partial}{\partial m}V(m)(x)\Delta m(x)dx+\int_{\mathbb{R}^{n}}\int_{\mathbb{R}^{n}}\frac{\partial^{2}}{\partial m^{2}}V(m)(x,y)Dm(x).Dm(y)dxdy
\]

which can be rewritten as 

\begin{equation}
\sum_{k=1}^{n}D\mathcal{U}(X)(e_{k},e_{k})=\int_{\mathbb{R}^{n}}U(x,m)\Delta m(x)dx+\int_{\mathbb{R}^{n}}\int_{\mathbb{R}^{n}}\frac{\partial}{\partial m}U(x,m)(y)Dm(x).Dm(y)dxdy\label{eq:5-63}
\end{equation}

\[
=\int_{\mathbb{R}^{n}}\Delta U(x,m)m(x)dx+\int_{\mathbb{R}^{n}}\int_{\mathbb{R}^{n}}\frac{\partial}{\partial m}U(x,m)(y)Dm(x).Dm(y)dxdy
\]
 If we differentiate in $X,$ the left hand side gives $\sum_{k=1}^{n}D^{2}\mathcal{U}(X)(e_{k},e_{k}),$
which we want to compute. The right hand side is a functional of $m.$
So we first take the Gâteaux differential in $m,$ which is 

\begin{align*}
\Delta U(x,m)+\int_{\mathbb{R}^{n}}\Delta_{\xi}\frac{\partial}{\partial m}U(\xi,m)(x)m(\xi)d\xi+
\end{align*}
\[
+\int_{\mathbb{R}^{n}}\int_{\mathbb{R}^{n}}\frac{\partial^{2}U(x,m,t)}{\partial m^{2}}(\xi,\eta)Dm(\xi)Dm(\eta)d\xi d\eta-2\text{div}\left(\int_{\mathbb{R}^{n}}\frac{\partial U(x,m,t)}{\partial m}(\eta)Dm(\eta)d\eta\right)
\]

which is also 

\begin{align*}
\Delta_{x}U(x,m)+\int_{\mathbb{R}^{n}}\frac{\partial}{\partial m}U(x,m)(\xi)\Delta m(\xi)d\xi+
\end{align*}
\[
+\int_{\mathbb{R}^{n}}\int_{\mathbb{R}^{n}}\frac{\partial^{2}U(x,m,t)}{\partial m^{2}}(\xi,\eta)Dm(\xi)Dm(\eta)d\xi d\eta-2\text{div}\left(\int_{\mathbb{R}^{n}}\frac{\partial U(x,m,t)}{\partial m}(\eta)Dm(\eta)d\eta\right)
\]
 We then have to take the gradient in $x,$ and replace $x$ by $X.$
We obtain immediately the relation (\ref{eq:5-62}).

\subsection{SYSTEM OF HJB-FP EQUATIONS }

As in the deterministic case $\beta=0,$ see section \ref{sub:INTERPRETATION-OF-THE}
we can derive from the Master equation a system of coupled stochastic
HJB-FP equations. They reduce to the deterministic system (\ref{eq:2.12})
to (\ref{eq:2.15}) when $\beta=0.$ Consider the conditional probability
density corresponding to the optimal feedback $\hat{v}(x,m,DU(x,m,t))$,
see (\ref{eq:5-1}), solution of 

\begin{equation}
\partial_{t}m+(A^{*}m-\frac{1}{2}\beta^{2}\Delta m+\text{div}(G(x,m,DU)m))dt+\beta Dm.db(t)=0\label{eq:3.6}
\end{equation}

\[
m(x,0)=m_{0}(x)
\]
 and set $u(x,t)=U(x,m(t),t)$, we obtain 

\begin{align}
-\partial_{t}u+(Au-\frac{1}{2}\beta^{2}\Delta u)dt & +\beta^{2}\text{div}\left(\int_{\mathbb{R}^{n}}\frac{\partial U(x,m,t)}{\partial m}(\xi)Dm(\xi)d\xi\right)dt=\label{eq:3.7}\\
H(x,m,Du(x))+ & \int_{\mathbb{R}^{n}}\frac{\partial}{\partial m}H(\xi,m,Du(\xi))(x)m(\xi)d\xi+\beta\int_{\mathbb{R}^{n}}\frac{\partial U(x,m,t)}{\partial m}(\xi)Dm(\xi)d\xi db(t)\nonumber \\
u(x,T) & =h(x,m)+\int_{\mathbb{R}^{n}}\frac{\partial}{\partial m}h(\xi,m)(x)m(\xi)d\xi\label{eq:3.8}
\end{align}

Setting 

\[
B(x,t)=\int_{\mathbb{R}^{n}}\frac{\partial U(x,m,t)}{\partial m}(\xi)Dm(\xi)d\xi
\]
 we can rewrite (\ref{eq:3.6}), (\ref{eq:3.7}) as follows 

\begin{align}
-\partial_{t}u+(Au-\frac{1}{2}\beta^{2}\Delta u)dt & +\beta^{2}\text{div}\: B(x,t)dt=\label{eq:3.7-1}\\
H(x,m,Du(x))+ & \int_{\mathbb{R}^{n}}\frac{\partial}{\partial m}H(\xi,m,Du(\xi))(x)m(\xi)d\xi+\beta B(x,t)db(t)\nonumber \\
u(x,T) & =h(x,m)+\int_{\mathbb{R}^{n}}\frac{\partial}{\partial m}h(\xi,m)(x)m(\xi)d\xi\label{eq:3.8-1}
\end{align}

\begin{equation}
\partial_{t}m+(A^{*}m-\frac{1}{2}\beta^{2}\Delta m+\text{div}(G(x,m,Du)m))dt+\beta Dm.db(t)=0\label{eq:3.6-1}
\end{equation}

\[
m(x,0)=m_{0}(x)
\]
 In fact the equation for $u$ contains two unknowns $u(x,t)$ and
$B(x,t)$ which are both field processes adapted to the filtration
$\mathcal{B}^{t}.$ Since the equation for $u$ is a backward stochastic
partial differential equation ( backward SPDE) we need the two unknowns
to obtain adapted processes.

\section{\label{sub:MEAN-FIELD-GAMES}MEAN FIELD GAMES }

\subsection{THE PROBLEM}

In Mean Field Games, we cannot have a Bellman equation, similar to
(\ref{eq:2.6}), (\ref{eq:2.7}), since the problem is not a control
problem. However, for a fixed parameter $m(.)$ we can introduce 

\begin{align}
dx & =g(x,m,v(x))ds+\sigma(x)dw\label{eq:1.3-1}\\
x(t) & =x\nonumber 
\end{align}

\begin{equation}
J_{x,t}(v(.),m(.))=\mathbb{E}[\int_{t}^{T}f(x(s),m(s),v(x(s)))ds+h(x(T),m(T))]\label{eq:1.4-1}
\end{equation}

and we set 

\[
u(x,t)=\inf_{v(.)}J_{x,t}(v(.),m(.))
\]
 in which we omit to write explicitly the dependence in $m.$ Then
$u(x,t)$ satisfies Bellman equation 

\begin{align}
-\frac{\partial u}{\partial t}+Au & =H(x,m,Du(x))\label{eq:2.12-1}\\
u(x,T) & =h(x,m)\label{eq:2.13-1}
\end{align}

We next write that $m$ must be the probability density of the optimal
state, hence 

\begin{align}
\frac{\partial m}{\partial t}+A^{*}m+\text{div }(G(x,m,Du)m(x)) & =0\label{eq:2.14-1}\\
m(x,0)=m_{0}(x)\label{eq:2.15-1}
\end{align}
 and this is the system of HJB-FP equations, corresponding to the
classical Mean Field Games problem.

\subsection{MASTER EQUATION}

We do not have Dynamic Programming, however the uncoupling argument
corresponding to the system (\ref{eq:2.12-1}),(\ref{eq:2.13-1}),(\ref{eq:2.14-1}),
(\ref{eq:2.15-1}) can be pursued. This time, we introduce directly
the Master equation by setting

\begin{align}
-\frac{\partial U}{\partial t}+AU+\int_{\mathbb{R}^{n}}\frac{\partial}{\partial m}U(x,m,t)(\xi)(A^{*}m(\xi)+\text{div (}G(\xi,m,DU(\xi))m(\xi))d\xi & =H(x,m,DU(x))\label{eq:2.10-1}
\end{align}
\begin{equation}
U(x,m,T)=h(x,m)\label{eq:2.11-1}
\end{equation}

We then consider $m(x,t)$ solution of the FP equation 

\begin{equation}
\frac{\partial m}{\partial t}+A^{*}m+\text{div }(G(x,m,DU(x,m))m(x))=0\label{eq:2.11-2}
\end{equation}
\[
m(x,0)=m_{0}(x)
\]

We next introduce

\begin{equation}
u(x,t)=U(x,m(t),t)\label{eq:2.12-2}
\end{equation}

in which $m(t)\equiv m(x,t)$ is the solution of (\ref{eq:2.11-2}).
In combining (\ref{eq:2.10-1}) and (\ref{eq:2.11-2}) we obtain easily
(\ref{eq:2.12-1}), in which $m(t)$ is the solution of (\ref{eq:2.11-2}).
The master equation (\ref{eq:2.10-1}) looks very similar to that
of mean field type control, see (\ref{eq:2.10}), (\ref{eq:2.11}),
and even simpler, since the derivative in $m$ of the Hamiltonian
does not appear anymore. However, we may loose a very important symmetry
property 

\begin{equation}
\frac{\partial}{\partial m}U(x,m,t)(\xi)\neq\frac{\partial}{\partial m}U(\xi,m,t)(x)\label{eq:2.122}
\end{equation}

\subsection{MASTER EQUATION IN $\mathcal{H}$ }

We want to associate to the Master equation (\ref{eq:2.10-1}) (\ref{eq:2.11-1})
an equation in the space $\mathcal{H}.$ We do not have a Bellman
equation, so we shall proceed directly with (\ref{eq:2.10-1}), (\ref{eq:2.11-1}).
To simplify we assume $\sigma(x)=I,$ so we rewrite the Master equation
as 

\begin{equation}
\frac{\partial U}{\partial t}+\frac{\sigma^{2}}{2}\Delta_{x}U+\frac{\sigma^{2}}{2}\int_{\mathbb{R}^{n}}\Delta_{\xi}\frac{\partial}{\partial m}U(x,m,t)(\xi)m(\xi)d\xi+\int_{\mathbb{R}^{n}}D_{\xi}\frac{\partial}{\partial m}U(x,m,t)(\xi).G(\xi,m,D_{\xi}U(\xi)))m(\xi)d\xi+\label{eq:5-6-2}
\end{equation}

\[
+H(x,m,D_{x}U(x))=0
\]
\[
U(x,m,T)=h(x,m)
\]
 We next take the gradient in $x$ to obtain 

\begin{equation}
\frac{\partial}{\partial t}D_{x}U+\frac{\sigma^{2}}{2}\Delta_{x}D_{x}U+\frac{\sigma^{2}}{2}D_{x}\int_{\mathbb{R}^{n}}\Delta_{\xi}\frac{\partial}{\partial m}U(x,m,t)(\xi)m(\xi)d\xi+D_{x}\int_{\mathbb{R}^{n}}D_{\xi}\frac{\partial}{\partial m}U(x,m,t)(\xi).G(\xi,m,DU(\xi)))m(\xi)d\xi+\label{eq:5-7-1}
\end{equation}
\[
+D_{x}H(x,m,D_{x}U(x))+D_{x}^{2}U(x,m,t)G(x,m,D_{x}U)=0
\]
 We then replace $x$ by $X$ and we want to interpret the equation
which is obtained. We set $\mathcal{U}(X,t)=D_{x}U(X,\mathcal{L}_{X},t)$,
then of course $\dfrac{\partial}{\partial t}D_{x}U(X,\mathcal{L}_{X},t)=\dfrac{\partial}{\partial t}\mathcal{U}$.
Next 

\[
D_{x}H(X,m,D_{x}U(X))=D_{x}H(X,\mathcal{L}_{X},\mathcal{U}(X))
\]
 We have next , see (\ref{eq:5-11})

\begin{align}
D_{x}\int_{\mathbb{R}^{n}}D_{\xi}\frac{\partial}{\partial m}U(X,m,t)(\xi).G(\xi,m,D_{\xi}U(\xi)))m(\xi)d\xi+D_{x}^{2}U(X,m,t)G(X,m,D_{x}U(X)) & =\label{eq:5-11-1}\\
D\mathcal{U}(X,t)G(X,\mathcal{L}_{X},\mathcal{U}(X))\nonumber 
\end{align}
 There remains to interpret $\Delta_{x}D_{x}U(X)+D_{x}\int_{\mathbb{R}^{n}}\Delta_{\xi}\frac{\partial}{\partial m}U(X,m,t)(\xi)m(\xi)d\xi$.
We first have 

\[
(D\mathcal{U}(X)\Gamma,\Gamma)=\int_{\mathbb{R}^{n}}\Delta_{x}U(X,m)m(x)dx
\]
\[
D_{X}((D\mathcal{U}(X)\Gamma,\Gamma)=(D^{2}\mathcal{U}(X)\Gamma,\Gamma)=
\]
\[
D_{x}\Delta_{x}U(X,\mathcal{L}_{X})+D_{x}\int_{\mathbb{R}^{n}}\Delta_{\xi}\frac{\partial}{\partial m}U(\xi,m,t)(X)m(\xi)d\xi
\]
which is the quantity we want only if we make the symmetry assumption 

\begin{equation}
\frac{\partial}{\partial m}U(x,m,t)(\xi)=\frac{\partial}{\partial m}U(\xi,m,t)(x)\label{eq:5-112}
\end{equation}

So under the assumption (\ref{eq:5-112}), we claim that $\mathcal{U}(X,t)$
is the solution of the Master equation 

\begin{equation}
\frac{\partial}{\partial t}\mathcal{U}+\frac{\sigma^{2}}{2}D^{2}\mathcal{U}(\Gamma,\Gamma)+D\mathcal{U}(X)G(X,\mathcal{L}_{X},\mathcal{U}(X))+D_{x}H(X,\mathcal{L}_{X},\mathcal{U}(X))=0\label{eq:5-61-2}
\end{equation}

\[
\mathcal{U}(X,T)=D_{x}h(X,\mathcal{L}_{X})
\]
 However, (\ref{eq:5-112}) implies that $D$$\mathcal{U}$ is self-adoint.
So we are looking to a solution $\mathcal{U}(X,t)$ such that $D\mathcal{U}(X,t)$
is a self-adjoint operator in $\mathcal{H}$ . As we shall see in
the linear quadratic case, this might not be true.

Note that , applying (\ref{eq:5-112}) at time $T$ we need to have,
in particular 

\begin{equation}
\frac{\partial}{\partial m}h(x,m)(\xi)=\frac{\partial}{\partial m}h(\xi,m)(x)\label{eq:5-613}
\end{equation}

In the mean-field control problem, since $\mathcal{U}(X,T)=D$$h(X,\mathcal{L}_{X})$,
we have $D\mathcal{U}(X,T)=D^{2}h(X,\mathcal{L}_{X}),$ it is a self
adjoint operator by construction. 

Therefore there remains an advantage for equation (\ref{eq:5-6-2}).
There is still a possibility to use the space $\mathcal{H},$ as follows.
Instead of considering the derivative $D_{x}U(x,m)$ one considers
the map $U(x,m)$ itself. We then write $U(x,m)=U(x,\mathcal{L}_{X}).$
One then interpret directly the Master equation (\ref{eq:5-6-2}).
We have successively 

\[
\int_{\mathbb{R}^{n}}\Delta_{\xi}\frac{\partial}{\partial m}U(x,m,t)(\xi)m(\xi)d\xi=D_{X}^{2}U(x,X,t)(\Gamma,\Gamma)
\]
\[
\int_{\mathbb{R}^{n}}D_{\xi}\frac{\partial}{\partial m}U(x,m,t)(\xi).G(\xi,m,D_{\xi}U(\xi)))m(\xi)d\xi=(D_{X}U(x,X,t),G(X,\mathcal{L}_{X},D_{x}U(X,\mathcal{L}_{X},t))
\]
 so we get the Master equation 

\begin{equation}
\frac{\partial U}{\partial t}+\frac{\sigma^{2}}{2}\Delta_{x}U+\frac{\sigma^{2}}{2}D_{X}^{2}U(x,X,t)(\Gamma,\Gamma)+(D_{X}U(x,X,t),G(X,\mathcal{L}_{X},D_{x}U(X,\mathcal{L}_{X},t))\label{eq:5-61-1-1}
\end{equation}

\[
+H(x,\mathcal{L}_{X},D_{x}U(X,\mathcal{L}_{X},t))=0
\]

\[
U(x,X,T)=h(X,\mathcal{L}_{X})
\]

\section{STOCHASTIC MEAN FIELD GAMES}

\subsection{THE MODEL }

We recall how to obtain a system of HJB-FP equations, see \cite{BFY1}
for more details. 

Consider the notation of section \ref{sub:PRELIMINARIES}. We recall
that the state equation is the solution of 

\begin{align}
dx & =g(x,m(t),v(x,t))dt+\sigma(x)dw+\beta db(t)\label{eq:4.1}\\
x(0) & =x_{0}\nonumber 
\end{align}

and $\mathcal{B}^{t}$=$\sigma(b(s),s\leq t)$ , $\mathcal{F}$$^{t}$=$\sigma(x_{0},b(s),w(s),\: s\leq t)$.
This time $m(t)$ is a given process adapted to $\mathcal{B}^{t}$
with values in $L^{2}(\mathbb{R}^{n}).$ We again consider feedback controls
which are field processes adapted to the filtration $\mathcal{B}^{t}$. 

We follow the theory developed by Shi Ge Peng \cite{SGP}. We define
the function 

\begin{equation}
u(x,t)=\inf_{v(.)}\: \mathbb{E}^{\mathcal{B}^{t}}[\int_{t}^{T}f(x(s),m(s),v(x(s),s))ds+h(x(T),m(T))]\label{eq:4.2}
\end{equation}
and show that it satisfies a stochastic Hamilton-Jacobi-Bellman equation.
The equation is 

\begin{align}
-\partial_{t}u+(Au-\frac{1}{2}\beta^{2}\Delta u)dt & +\beta^{2}\text{div}\: B(x,t)dt=H(x,m(t),Du(x))dt+\beta B(x,t)db(t)\label{eq:4.3}\\
u(x,T) & =h(x,m(T))\label{eq:4.4}
\end{align}
If we can solve (\ref{eq:4.3}), (\ref{eq:4.4}) for the pair $u(x,t),\: B(x,t)$
of adapted to $\mathcal{B}^{t}$ field processes, then a verification
argument has been described in \cite{BFY1}, to show that $u(x,t)$
coincides with (\ref{eq:4.2}). 

Next, consider the optimal feedback $\hat{v}(x,m(t),Du(x,t))$ and
impose the fixed point property that $m(t)$ is conditional probability
density of the optimal state, we get the stochastic FP equation 

\begin{equation}
\partial_{t}m+(A^{*}m-\frac{1}{2}\beta^{2}\Delta m+\text{div}(G(x,m,Du)m))dt+\beta Dm.db(t)=0\label{eq:4.5}
\end{equation}

\[
m(x,0)=m_{0}(x)
\]

We thus have obtained the pair of HJB-FP equations for the stochastic
mean field game problem, (\ref{eq:4.3}),(\ref{eq:4.4}) and (\ref{eq:4.5}).

\subsection{THE MASTER EQUATION}

We  derive the Master equation by writing 

\[
u(x,t)=U(x,m(t),t)
\]

The calculation has been done in \cite{BFY1}. We obtain

\begin{equation}
-\frac{\partial U}{\partial t}+AU-\frac{1}{2}\beta^{2}\Delta U+\int_{\mathbb{R}^{n}}\frac{\partial U(x,m,t)}{\partial m}(\xi)(A^{*}m(\xi)-\frac{1}{2}\beta^{2}\Delta m(\xi)+\text{div }(G(\xi,m,DU(\xi))m(\xi)))d\xi+\label{eq:4.7}
\end{equation}
\[
+\beta^{2}\text{div}\:(\int_{\mathbb{R}^{n}}\frac{\partial U(x,m,t)}{\partial m}(\xi)Dm(\xi)d\xi)-\frac{1}{2}\beta^{2}\int_{\mathbb{R}^{n}}\int_{\mathbb{R}^{n}}\frac{\partial^{2}U(x,m,t)}{\partial m^{2}}(\xi,\eta)Dm(\xi)Dm(\eta)d\xi d\eta=H(x,m,DU(x))
\]
\[
U(x,m,T)=h(x,m)
\]

\subsection{THE MASTER EQUATION IN THE SPACE $\mathcal{H}$ . }

We take the particular case $\sigma(x)=I$ then (\ref{eq:4.7}) becomes 

\begin{equation}
-\frac{\partial U}{\partial t}-\frac{1}{2}(\sigma^{2}+\beta^{2})\Delta_{x}U-\frac{1}{2}(\sigma^{2}+\beta^{2})\int_{\mathbb{R}^{n}}\Delta_{\xi}\frac{\partial U(x,m,t)}{\partial m}(\xi)m(\xi)d\xi-\int_{\mathbb{R}^{n}}D_{\xi}\frac{\partial U(x,m,t)}{\partial m}(\xi).G(\xi,m,DU(\xi))m(\xi)d\xi\label{eq:4-7}
\end{equation}

\[
+\beta^{2}\text{div}\:(\int_{\mathbb{R}^{n}}\frac{\partial U(x,m,t)}{\partial m}(\xi)Dm(\xi)d\xi)-\frac{1}{2}\beta^{2}\int_{\mathbb{R}^{n}}\int_{\mathbb{R}^{n}}\frac{\partial^{2}U(x,m,t)}{\partial m^{2}}(\xi,\eta)Dm(\xi)Dm(\eta)d\xi d\eta=H(x,m,DU(x))
\]

We assume the symmetry property (\ref{eq:5-112}). We take the gradient
in $x$ and obtain 

\begin{equation}
\frac{\partial}{\partial t}D_{x}U+\frac{\sigma^{2}+\beta^{2}}{2}\Delta_{x}D_{x}U+\frac{\sigma^{2}+\beta^{2}}{2}D_{x}\int_{\mathbb{R}^{n}}\Delta_{\xi}\frac{\partial}{\partial m}U(x,m,t)(\xi)m(\xi)d\xi+\label{eq:5-7-1-1}
\end{equation}
\[
+D_{x}\int_{\mathbb{R}^{n}}D_{\xi}\frac{\partial}{\partial m}U(x,m,t)(\xi).G(\xi,m,DU(\xi))m(\xi)d\xi+D_{x}H(x,m,D_{x}U(x))+D_{x}^{2}U(x,m,t)G(x,m,D_{x}U)-
\]
\[
-\beta^{2}\text{div}\:(\int_{\mathbb{R}^{n}}\frac{\partial U(x,m,t)}{\partial m}(\xi)Dm(\xi)d\xi)+\frac{1}{2}\beta^{2}\int_{\mathbb{R}^{n}}\int_{\mathbb{R}^{n}}\frac{\partial^{2}U(x,m,t)}{\partial m^{2}}(\xi,\eta)Dm(\xi)Dm(\eta)d\xi d\eta=0
\]

So introducing $\mathcal{U}(X)=D_{x}U(X,\mathcal{L}_{X})$ and we
replace $x$ by $X$ in the equation (\ref{eq:5-7-1-1}). We have
interpreted in (\ref{eq:5-61-2}) all the terms not including $\beta.$We
want to check that 

\begin{equation}
\Delta_{x}D_{x}U(X)+D_{x}\int_{\mathbb{R}^{n}}\Delta_{\xi}\frac{\partial}{\partial m}U(X,m,t)(\xi)m(\xi)d\xi-\label{eq:5-71}
\end{equation}

\[
-2D_{x}\text{div}\:(\int_{\mathbb{R}^{n}}\frac{\partial U(X,m,t)}{\partial m}(\xi)Dm(\xi)d\xi)+D_{x}\int_{\mathbb{R}^{n}}\int_{\mathbb{R}^{n}}\frac{\partial^{2}U(X,m,t)}{\partial m^{2}}(\xi,\eta)Dm(\xi)Dm(\eta)d\xi d\eta=\sum_{k=1}^{n}D^{2}\mathcal{U}(X)(e_{k},e_{k})
\]
 Now we have 

\begin{equation}
\sum_{k=1}^{n}D\mathcal{U}(X)(e_{k},e_{k})=\int_{\mathbb{R}^{n}}\Delta U(x,m)m(x)dx+\int_{\mathbb{R}^{n}}\int_{\mathbb{R}^{n}}\frac{\partial}{\partial m}U(x,m)(y)Dm(x).Dm(y)dxdy\label{eq:5-63-1}
\end{equation}

Next we take the derivative in $m$ of the right hand side and then
the gradient in $x$ to obtain

\begin{align*}
D_{x}\Delta_{x}U(x,m)+D_{x}\int_{\mathbb{R}^{n}}\Delta_{\xi}\frac{\partial}{\partial m}U(\xi,m,t)(x)m(\xi)d\xi+
\end{align*}
\[
+D_{x}\int_{\mathbb{R}^{n}}\int_{\mathbb{R}^{n}}\frac{\partial^{2}U(x,m,t)}{\partial m^{2}}(\xi,\eta)Dm(\xi)Dm(\eta)d\xi d\eta-2D_{x}\text{div}\left(\int_{\mathbb{R}^{n}}\frac{\partial U(x,m,t)}{\partial m}(\eta)Dm(\eta)d\eta\right)
\]

We need again to use the symmetry property and replace $x$ by $X$
to prove (\ref{eq:5-71}). We can now write the Master equation in
the space $\mathcal{H}$ , namely 

\begin{equation}
\frac{\partial}{\partial t}\mathcal{U}+\frac{\sigma^{2}}{2}D^{2}\mathcal{U}(\Gamma,\Gamma)+\frac{\beta^{2}}{2}\sum_{k=1}^{n}D^{2}\mathcal{U}(X)(e_{k},e_{k})+D\mathcal{U}(X)G(X,\mathcal{L}_{X},\mathcal{U}(X))+D_{x}H(X,\mathcal{L}_{X},\mathcal{U}(X))=0\label{eq:5-61-2-1}
\end{equation}
\[
\mathcal{U}(X,T)=D_{x}h(X,\mathcal{L}_{X})
\]

If we do not have the symmetry property, then we consider the Master
equation for $U(x,m,t)$, namely (\ref{eq:4-7}), and we introduce
the argument $X$ instead of $m,$ as in (\ref{eq:5-61-1-1}). We
use (\ref{eq:1-40}) to state 

\[
\int_{\mathbb{R}^{n}}\int_{\mathbb{R}^{n}}\frac{\partial^{2}U(x,m,t)}{\partial m^{2}}(\xi,\eta)Dm(\xi)Dm(\eta)d\xi d\eta=\sum_{k=1}^{n}D_{X}^{2}U(x,X,t)(e_{k},e_{k})-D_{X}^{2}U(x,X,t)(\Gamma,\Gamma)
\]
 Next 

\[
\text{div}\:(\int_{\mathbb{R}^{n}}\frac{\partial U(x,m,t)}{\partial m}(\xi)Dm(\xi)d\xi)=-\sum_{k=1}^{n}\frac{\partial}{\partial x_{k}}(e_{k},D_{X}U(x,X,t))
\]
 We can then write the Master equation for $U(x,X,t)$

\begin{equation}
\frac{\partial U}{\partial t}+\frac{\sigma^{2}+\beta^{2}}{2}\Delta_{x}U+\frac{\sigma^{2}}{2}D_{X}^{2}U(x,X,t)\Gamma,\Gamma)+(D_{X}U(x,X,t),G(X,\mathcal{L}_{X},D_{x}U(X,\mathcal{L}_{X},t))+\label{eq:5-614-1}
\end{equation}

\[
+\frac{\beta^{2}}{2}\sum_{k=1}^{n}D_{X}^{2}U(x,X,t)(e_{k},e_{k})+\beta^{2}\sum_{k=1}^{n}\frac{\partial}{\partial x_{k}}(e_{k},D_{X}U(x,X,t))+H(x,\mathcal{L}_{X},D_{x}U(x,X,t))=0
\]

\[
U(x,X,T)=h(x,\mathcal{L}_{X})
\]

\section{LINEAR QUADRATIC PROBLEMS}

\subsection{ASSUMPTIONS AND GENERAL COMMENTS}

For linear quadratic problems, we know that explicit formulas can
be obtained. We have solved the Master equations in $U(x,m,t)$ in
our paper \cite{BFY1}. So we shall here solve only the Master equation
in the space $\mathcal{H}.$ We then take 

\begin{equation}
f(x,m,v)=\frac{1}{2}[x^{*}Qx+v^{*}Rv+(x-Sy)^{*}\bar{Q}(x-Sy)]\label{eq:5.2}
\end{equation}

\begin{equation}
g(x,m,v)=Ax+\bar{A}y+Bv\label{eq:5.3}
\end{equation}

\begin{equation}
h(x,m)=\frac{1}{2}[x^{*}Q_{T}x+(x-S_{T}y)^{*}\bar{Q}_{T}(x-S_{T}y)]\label{eq:5.4}
\end{equation}

in which we have noted $y=\int_{\mathbb{R}^{n}}\xi m(\xi)d\xi$ . We also assume 

\begin{equation}
\sigma(x)=\sigma,\;\text{hence}\; a(x)=a=\sigma\sigma^{*}\label{eq:5.5}
\end{equation}

We deduce easily 

\begin{align}
H(x,m,q) & =\frac{1}{2}x^{*}(Q+\bar{Q})\, x-x^{*}\bar{Q}Sy+\frac{1}{2}y^{*}S^{*}\bar{Q}Sy\label{eq:5.6}\\
- & \frac{1}{2}q^{*}BR^{-1}B^{*}q+q^{*}(Ax+\bar{A}y)\nonumber 
\end{align}

\begin{equation}
G(x,m,q)=Ax+\bar{A}y-BR^{-1}B^{*}q\label{eq:5.7}
\end{equation}

\subsection{MEAN FIELD TYPE CONTROL MASTER EQUATION}

We begin with Bellman equation (\ref{eq:2-8}), namely

\begin{equation}
\frac{\partial V}{\partial t}+\frac{1}{2}D^{2}V(X)(\sigma(X)\Gamma,\sigma(X)\Gamma)+\frac{1}{2}\beta^{2}\sum_{k=1}^{n}D^{2}V(X)(e_{k},e_{k})+\mathbb{E}H(X,\mathcal{L}_{X},DV(X))=0\label{eq:2-8-1}
\end{equation}
\[
V(X,T)=\mathbb{E}h(X,\mathcal{L}_{X})
\]
 and using (\ref{eq:5-6}) we get 

\begin{equation}
\frac{\partial V}{\partial t}+\frac{1}{2}D^{2}V(X)(\sigma(X)\Gamma,\sigma(X)\Gamma)+\frac{1}{2}\beta^{2}\sum_{k=1}^{n}D^{2}V(X)(e_{k},e_{k})+\label{eq:2-81}
\end{equation}

\[
+\mathbb{E}\left(\frac{1}{2}X^{*}(Q+\bar{Q})\, X-X^{*}\bar{Q}S\, \mathbb{E}X+\frac{1}{2}\mathbb{E}X^{*}S^{*}\bar{Q}S\, \mathbb{E}X-\frac{1}{2}DV^{*}BR^{-1}B^{*}DV+DV^{*}(AX+\bar{A}\mathbb{E}X)\right)=0
\]
\[
V(X,T)=\frac{1}{2}\mathbb{E}\,[X^{*}Q_{T}X+(X-S_{T}\mathbb{E}X)^{*}\bar{Q}_{T}(X-S_{T}\mathbb{E}X)]
\]

We look for a solution of (\ref{eq:2-81}), of the form 

\begin{equation}
V(X,t)=\frac{1}{2}\mathbb{E}\, X^{*}P(t)X+\frac{1}{2}(\mathbb{E}X)^{*}\Sigma(t)\mathbb{E}X+\lambda(t)\label{eq:5.13}
\end{equation}

We have clearly 

\begin{align}
P(T) & =Q_{T}+\bar{Q_{T}},\;\Sigma(T)=S_{T}^{*}\bar{Q}_{T}S_{T}-(S_{T}^{*}\bar{Q}_{T}+\bar{Q}_{T}S_{T}),\:\lambda(T)=0\label{eq:5.14}
\end{align}

Next 

\[
DV(X,t)=P(t)X+\Sigma(t)\, \mathbb{E}X
\]
\[
D^{2}V(X,t)Z=P(t)Z+\Sigma(t)\, \mathbb{E}Z
\]

therefore 

\[
D^{2}V(X)(\sigma(X)\Gamma,\sigma(X)\Gamma)=\text{tr}P(t)a
\]

\[
\sum_{k=1}^{n}D^{2}V(X)(e_{k},e_{k})=\text{tr }(P(t)+\Sigma(t))
\]

With these calculations, we can proceed in equation (\ref{eq:2-81})
and obtain 

\[
\frac{1}{2}\mathbb{E}\, X^{*}\frac{d}{dt}P(t)X+\frac{1}{2}(\mathbb{E}X)^{*}\frac{d}{dt}\Sigma(t)\mathbb{E}X+\frac{d}{dt}\lambda(t)
\]
\[
+\frac{1}{2}\text{tr}\, aP(t)+\frac{\beta^{2}}{2}\text{tr}(P(t)+\Sigma(t))+
\]
\[
+\mathbb{E}\left[\frac{1}{2}X^{*}(Q+\bar{Q})X-X^{*}\bar{Q}S\, \mathbb{E}X+\frac{1}{2}(\mathbb{E}X)^{*}S^{*}\bar{Q}S\mathbb{E}X-\right.
\]
\[
\left.-\frac{1}{2}(P(t)X+\Sigma(t)\mathbb{E}X)^{*}BR^{-1}B^{*}(P(t)X+\Sigma(t)\mathbb{E}X)+(P(t)X+\Sigma(t)\mathbb{E}X)^{*}(AX+\bar{A}\mathbb{E}X)\right]=0
\]
 We can identify terms and obtain 

\[
\frac{d\lambda(t)}{dt}+\frac{1}{2}\text{tr}\, aP(t)+\frac{\beta^{2}}{2}\text{tr}(P(t)+\Sigma(t))=0
\]
 therefore, from the final condition $\lambda(T)=0,$ it follows

\begin{equation}
\lambda(t)=\int_{t}^{T}(\frac{1}{2}\text{tr}\, aP(s)+\frac{\beta^{2}}{2}\text{tr}(P(s)+\Sigma(s))ds\label{eq:5.16}
\end{equation}
Identifying quadratic terms in $X$( within the expected value) and
in $\mathbb{E}X$ respectively, it follows easily that 

\begin{equation}
\frac{dP}{dt}+PA+A^{*}P-PBR^{-1}B^{*}P+Q+\bar{Q}=0,\; P(T)=Q_{T}+\bar{Q_{T}},\label{eq:5.17}
\end{equation}
 and 

\begin{equation}
\frac{d\Sigma}{dt}+\Sigma(A+\bar{A}-BR^{-1}B^{*}P)+(A+\bar{A}-BR^{-1}B^{*}P)^{*}\Sigma-\label{eq:5.18}
\end{equation}

\[
-\Sigma BR^{-1}B^{*}\Sigma+S^{*}\bar{Q}S-\bar{Q}S-S^{*}\bar{Q}+P\bar{A}+\bar{A}\,^{*}P=0
\]
\[
\Sigma(T)=S_{T}^{*}\bar{Q}_{T}S_{T}\,-(S_{T}^{*}\bar{Q}_{T}+\bar{Q}_{T}S_{T})
\]

We obtain formula (\ref{eq:5.13}) with the values of $P(t),\,\Sigma(t),\,\lambda(t)$
given by equations (\ref{eq:5.17}), (\ref{eq:5.18}), (\ref{eq:5.16}). 

We turn to the Master equation, namely (\ref{eq:5-61-1}) 

\begin{equation}
\frac{\partial}{\partial t}\mathcal{U}+\frac{\sigma^{2}}{2}D^{2}\mathcal{U}(\Gamma,\Gamma)+\frac{1}{2}\beta^{2}\sum_{k=1}^{n}D^{2}\mathcal{U}(X)(e_{k},e_{k})+D\mathcal{U}(X)G(X,\mathcal{L}_{X},\mathcal{U}(X))+\label{eq:5-61-1-2}
\end{equation}

\[
D_{x}H(X,\mathcal{L}_{X},\mathcal{U}(X))+\mathbb{E}_{Y,\mathcal{U}(Y)}D_{x}\frac{\partial H}{\partial m}(Y,\mathcal{L}_{X},\mathcal{U}(Y))(X)=0
\]

\[
\mathcal{U}(X,T)=Dh(X,\mathcal{L}_{X})
\]
 We have 

\[
G(X,\mathcal{L}_{X},\mathcal{U}(X))=AX+\bar{A}\mathbb{E}X-BR^{-1}B^{*}\mathcal{U}(X)
\]
\[
D_{x}H(X,\mathcal{L}_{X},\mathcal{U}(X))=(Q+\bar{Q})X-\bar{Q}S\, \mathbb{E}X\,+A^{*}\mathcal{U}(X)
\]

\[
\frac{\partial H}{\partial m}(x,m,q)(\xi)=-x^{*}\bar{Q}S\xi+(\int\eta m(\eta)d\eta)^{*}S^{*}\bar{Q}S\xi+q^{*}\bar{A}\xi
\]
 so
\[
\mathbb{E}_{Y,\mathcal{U}(Y)}D_{x}\frac{\partial H}{\partial m}(Y,\mathcal{L}_{X},\mathcal{U}(Y))(X)=(-S^{*}\bar{Q}+S^{*}\bar{Q}S)\mathbb{E}X+(\bar{A})^{*}\mathbb{E}\mathcal{U}(X)
\]
 and (\ref{eq:5-61-1-2}) reads 

\begin{equation}
\frac{\partial}{\partial t}\mathcal{U}+\frac{\sigma^{2}}{2}D^{2}\mathcal{U}(\Gamma,\Gamma)+\frac{1}{2}\beta^{2}\sum_{k=1}^{n}D^{2}\mathcal{U}(X)(e_{k},e_{k})+(AX+\bar{A}\mathbb{E}X-BR^{-1}B^{*}\mathcal{U}(X))\label{eq:5-62-1}
\end{equation}
\[
+(Q+\bar{Q})X+(-\bar{Q}S-S^{*}\bar{Q}+S^{*}\bar{Q}S)\, \mathbb{E}X\,+A^{*}\mathcal{U}(X)+(\bar{A})^{*}\mathbb{E}\mathcal{U}(X)=0
\]
 We expect $\mathcal{U}(X,t)=P(t)X+\Sigma(t)\, \mathbb{E}X$ to be the solution.
This is satisfied at time $T.$ We check easily that equation (\ref{eq:5-62-1})
is satisfied with the choices of $P(t)$ and $\Sigma(t)$ given by
(\ref{eq:5.17}), (\ref{eq:5.18}), noting that $D\mathcal{U}(X,t)Z=P(t)Z+\Sigma(t)\mathbb{E}Z$.
We see also that (\ref{eq:5-62-1}) could be consider directly, without
referring to Bellman equation (\ref{eq:2-81}).

\subsection{MEAN FIELD GAMES MASTER EQUATION }

This time, there is no Bellman equation. We may consider the Master
equation 

\begin{equation}
\frac{\partial}{\partial t}\mathcal{U}+\frac{\sigma^{2}}{2}D^{2}\mathcal{U}(\Gamma,\Gamma)+\frac{\beta^{2}}{2}\sum_{k=1}^{n}D^{2}\mathcal{U}(X)(e_{k},e_{k})+D\mathcal{U}(X)G(X,\mathcal{L}_{X},\mathcal{U}(X))+D_{x}H(X,\mathcal{L}_{X},\mathcal{U}(X))=0\label{eq:5-62-2}
\end{equation}
\[
\mathcal{U}(X,T)=D_{x}h(X,\mathcal{L}_{X})
\]

which becomes in the LQ case 

\begin{equation}
\frac{\partial}{\partial t}\mathcal{U}+\frac{\sigma^{2}}{2}D^{2}\mathcal{U}(\Gamma,\Gamma)+\frac{\beta^{2}}{2}\sum_{k=1}^{n}D^{2}\mathcal{U}(X)(e_{k},e_{k})+D\mathcal{U}(X)(AX+\bar{A}\mathbb{E}X-BR^{-1}B^{*}\mathcal{U}(X))+\label{eq:5-63-3}
\end{equation}

\[
+(Q+\bar{Q})X-\bar{Q}S\, \mathbb{E}X\,+A^{*}\mathcal{U}(X)=0
\]
\[
\mathcal{U}(X,T)=(Q_{T}+\bar{Q_{T}})X-\bar{Q_{T}}S_{T}\mathbb{E}X
\]
 But we must have $D\mathcal{U}(X,t)$ self adjoint. At time $T$
we have 

\[
D\mathcal{U}(X,T)Z=(Q_{T}+\bar{Q_{T}})Z-\bar{Q_{T}}S_{T}\mathbb{E}Z
\]
 and therefore we need to assume 

\begin{equation}
\bar{Q_{T}}S_{T}=S_{T}^{*}\bar{Q_{T}}\label{eq:5-63-2}
\end{equation}

Next if we look for a linear solution $\mathcal{U}(X,t)=P(t)X+\Sigma(t)\mathbb{E}X,$
then we must have $P(t)$ and $\Sigma(t)$ symmetric. Checking in
(\ref{eq:5-63-3}) we see that $P(t)$ satisfies the Riccati equation
(\ref{eq:5.17}), but $\Sigma(t)$ satisfies 

\begin{equation}
\frac{d\Sigma}{dt}+\Sigma(A+\bar{A}-BR^{-1}B^{*}P)+(A-BR^{-1}B^{*}P)^{*}\Sigma-\label{eq:5.18-1}
\end{equation}

\[
-\Sigma BR^{-1}B^{*}\Sigma-\bar{Q}S+P\bar{A}=0
\]
\[
\Sigma(T)=-\bar{Q}_{T}S_{T}
\]
 and therefore $\Sigma(t)$ cannot be symmetric unless 

\begin{equation}
\bar{A}=0,\:\bar{Q}S=S^{*}\bar{Q}\label{eq:5-65}
\end{equation}

If these assumptions are not satisfied the equation (\ref{eq:5-63-3})
has no solution satisfying $D\mathcal{U}(X,t)$ self adjoint. In that
case we may use (\ref{eq:5-614-1}) which in the LQ case reduces to 

\begin{equation}
\frac{\partial U}{\partial t}+\frac{\sigma^{2}+\beta^{2}}{2}\Delta_{x}U+\frac{\sigma^{2}}{2}D_{X}^{2}U(x,X,t)(\Gamma,\Gamma)+(D_{X}U(x,X,t),AX+\bar{A}\mathbb{E}X-BR^{-1}B^{*}D_{x}U(X,\mathcal{L}_{X},t))+\label{eq:5-614-1-1}
\end{equation}

\[
+\frac{\beta^{2}}{2}\sum_{k=1}^{n}D_{X}^{2}U(x,X,t)(e_{k},e_{k})+\beta^{2}\sum_{k=1}^{n}\frac{\partial}{\partial x_{k}}(e_{k},D_{X}U(x,X,t))+\frac{1}{2}x^{*}(Q+\bar{Q})\, x-x^{*}\bar{Q}S\, \mathbb{E}X+
\]

\[
+\frac{1}{2}\mathbb{E}X^{*}S^{*}\bar{Q}S\, \mathbb{E}X-\frac{1}{2}(D_{x}U(x,X,t))^{*}BR^{-1}B^{*}D_{x}U(x,X,t)+(D_{x}U(x,X,t))^{*}(AX+\bar{A}\mathbb{E}X)=0
\]

\[
U(x,X,T)=\frac{1}{2}x^{*}(Q_{T}+\bar{Q}_{T})x-x^{*}\bar{Q}_{T}S_{T}\mathbb{E}X+\frac{1}{2}\mathbb{E}X^{*}S_{T}^{*}\bar{Q}_{T}S_{T}\mathbb{E}X
\]

We look for a solution of the form 

\begin{equation}
U(x,X,t)=\frac{1}{2}x^{*}P(t)x+x^{*}\Sigma(t)\mathbb{E}X+\frac{1}{2}\mathbb{E}X^{*}\Gamma(t)\mathbb{E}X+\mu(t)\label{eq:5.33}
\end{equation}

We must be careful that we cannot have $\Sigma(t)$ symmetric. It
is already true at time $T.$ We have 

\begin{equation}
P(T)=Q_{T}+\bar{Q}_{T}\label{eq:5.34}
\end{equation}

\begin{equation}
\Sigma(T)=-\bar{Q}_{T}S_{T}\label{eq:5.35}
\end{equation}
\begin{equation}
\Gamma(T)=S_{T}^{*}\bar{Q}_{T}S_{T}\label{eq:5.36}
\end{equation}

\begin{equation}
\mu(T)=0\label{eq:5.37}
\end{equation}
 Next

\[
D_{x}U(x,X,t)=P(t)x+\Sigma(t)\mathbb{E}X
\]
\[
D_{X}U(x,X,t)=\Sigma^{*}(t)x+\Gamma(t)\mathbb{E}X
\]
\[
D_{x}^{2}U(x,X,t)=P(t),\: D_{X}^{2}U(x,X,t)Z=\Gamma(t)\mathbb{E}Z
\]
 We use these formulas in the Master equation (\ref{eq:5-614-1-1}),
to obtain

\[
\frac{1}{2}x^{*}\frac{d}{dt}P(t)x+x^{*}\frac{d}{dt}\Sigma(t)\mathbb{E}X+\frac{1}{2}(\mathbb{E}X)^{*}\frac{d}{dt}\Gamma(t)\mathbb{E}X+\frac{d}{dt}\mu(t)+\frac{1}{2}(\beta^{2}+\sigma^{2})\text{tr}\, P(t)
\]

\[
+(\Sigma^{*}(t)x+\Gamma(t)\mathbb{E}X)^{*}(A+\bar{A}-BR^{-1}B^{*}(P(t)+\Sigma(t))\mathbb{E}X)+\frac{1}{2}\beta^{2}\text{tr}\Gamma(t)+\beta^{2}\text{tr}\Sigma(t)+
\]
\[
\frac{1}{2}x^{*}(Q+\bar{Q})x-x^{*}\bar{Q}S\mathbb{E}X+\frac{1}{2}(\mathbb{E}X)^{*}S^{*}\bar{Q}S(\mathbb{E}X)-
\]

\[
-\frac{1}{2}(P(t)x+\Sigma(t)\mathbb{E}X)^{*}BR^{-1}B^{*}(P(t)x+\Sigma(t)\mathbb{E}X)+(P(t)x+\Sigma(t)\mathbb{E}X)^{*}(Ax+\bar{A}\mathbb{E}X)=0
\]
 Identifying terms, we obtain 

\begin{equation}
\frac{d}{dt}P(t)+PA+A^{*}P-PBR^{-1}B^{*}P+Q+\bar{Q}=0\label{eq:5.38}
\end{equation}
\begin{equation}
\frac{d\Sigma}{dt}+\Sigma(A+\bar{A}-BR^{-1}B^{*}P)+(A^{*}-PBR^{-1}B^{*})\Sigma-\Sigma BR^{-1}B^{*}\Sigma-\bar{Q}S+P\bar{A}=0\label{eq:5.39}
\end{equation}
\begin{equation}
\frac{d\Gamma}{dt}+\Gamma(A+\bar{A}-BR^{-1}B^{*}(P+\Sigma))+(A+\bar{A}-BR^{-1}B^{*}(P+\Sigma))^{*}\Gamma+\label{eq:5.40}
\end{equation}
\[
+S^{*}\bar{Q}S-\Sigma BR^{-1}B^{*}\Sigma+\Sigma\bar{A}+\bar{A}\,^{*}\Sigma=0
\]
\begin{equation}
\frac{d}{dt}\mu(t)+\frac{\beta^{2}+\sigma^{2}}{2}\text{tr}P(t)+\frac{\beta^{2}}{2}\text{tr }\Gamma(t)+\beta^{2}\text{tr }\Sigma(t)=0\label{eq:5.41}
\end{equation}
 
\begin{rem}
\label{Rem3} In our previous work \cite{BFY1} we have solved completely
the Master equation in the LQ case, for mean field type control and
mean field games, considering $U(x,m,t)$ with $m\in L^{2}(\mathbb{R}^{n}).$
So $m$ was not necessarily the density of a probability. Although
everything can be carried out, the calculations are much more complex
than in the current case working with $U(x,X,t).$In particular we
had to keep a lot of terms depending on $m_{1}=\int_{\mathbb{R}^{n}}m(\xi)d\xi.$
In the case $U(x,X,t)$ this term is $1$ and calculations simplify
greatly. This shows the advantage of $\mathcal{H}=L^{2}(\Omega,\mathcal{A},P;\mathbb{R}^{n})$
with respect to $L^{2}(\mathbb{R}^{n}).$ We work with a Hilbert space in both
cases, but we keep the properties of probabilities in the case of
$\mathcal{H}.$ 
\end{rem}

\section{NEW CONTROL PROBLEMS}

\subsection{INTERPRETATION OF THE FIRST ORDER BELLMAN EQUATION }

We call first order Bellman equation, the equation 

\begin{align}
\frac{\partial V}{\partial t}+\mathbb{E}H(X,\mathcal{L}_{X},DV(X)) & =0\label{eq:9-1}\\
V(X,T)=\mathbb{E}\, h(X,\mathcal{L}_{X})\nonumber 
\end{align}
We associate to this equation a ``deterministic'' control problem
in the state space $\mathcal{H}.$ The state equation is given by 

\begin{align}
\frac{dX}{ds} & =g(X(s),\mathcal{L}_{X(s)},v(X(s),s))\label{eq:9-2}\\
X(t) & =X\nonumber 
\end{align}

The control is defined by a feedback $v(X,s)$ with values in $\mathbb{R}^{d}.$We
define the cost functional 

\begin{equation}
J_{X,t}(v(.))=\int_{t}^{T}\mathbb{E}f(X(s),\mathcal{L}_{X(s)},v(X(s),s))ds+\mathbb{E}h(X(T),\mathcal{L}_{X(T)})\label{eq:9-3}
\end{equation}

and the value function 

\begin{equation}
V(X,t)=\inf_{v(.)}\, J_{X,t}(v(.))\label{eq:9-4}
\end{equation}

We claim that the value function satisfies (\ref{eq:9-1}). Note that
although the equation (\ref{eq:9-2}) is deterministic, the trajectory
is random, because the initial condition $X$ is random. Nevertheless,
because the equation is deterministic, we say that (\ref{eq:9-2}),
( \ref{eq:9-3}) is a deterministic control problem in $\mathcal{H}$.
We provide a formal proof based on the optimality principle. We write 

\[
V(X,t)=\inf_{v(.)}\left[\int_{t}^{t+\epsilon}\mathbb{E}f(X(s),\mathcal{L}_{X(s)},v(X(s),s))ds+V(X(t+\epsilon),t+\epsilon)\right]
\]

\[
\sim\inf_{v}\left[\epsilon \mathbb{E}f(X,\mathcal{L}_{X},v)+V(X+\epsilon g(X,\mathcal{L}_{X},v),t+\epsilon)\right]
\]

from which we get easily 

\[
\frac{\partial V}{\partial t}+\inf_{v}\left[\mathbb{E}f(X,\mathcal{L}_{X},v)+(DV(X,t),g(X,\mathcal{L}_{X},v))\right]=0
\]

but 

\[
\inf_{v}\left[\mathbb{E}f(X,\mathcal{L}_{X},v)+(DV(X,t),g(X,\mathcal{L}_{X},v))\right]=\mathbb{E}\inf_{v}\left[f(X,\mathcal{L}_{X},v)+DV(X,t).g(X,\mathcal{L}_{X},v)\right]=
\]

\[
=\mathbb{E}H(X,\mathcal{L}_{X},DV(X,t))
\]
 and (\ref{eq:9-1}) follows immediately.

\subsection{INTERPRETATION OF THE FIRST STOCHASTIC BELLMAN EQUATION}

The first stochastic Bellman equation is 

\begin{align}
\frac{\partial V}{\partial t}+\frac{\sigma^{2}}{2}D^{2}V(X)(N,N)+\mathbb{E}\, H(X,\mathcal{L}_{X},\mathcal{U}(X)) & =0\label{eq:9-5}\\
V(X,T)=\mathbb{E}\, h(X,\mathcal{L}_{X})\nonumber 
\end{align}
where $N$ is standard gaussian independent of $X.$ We associate
to this equation a stochastic control problem in the space $\mathcal{H}$
. We assume that on the probability space $\Omega,\mathcal{A},P$
we can construct a standard Wiener process with values in $\mathbb{R}^{n}$
$w(t)$ which is independent of the random variable $X.$ Considering
again a feedback $v(X,s)$ we define the Ito equation in the space
$\mathcal{H}$ 

\begin{equation}
dX=g(X(s),\mathcal{L}_{X(s)},v(X(s),s))ds+\sigma dw(s)\label{eq:9-6}
\end{equation}
\[
X(t)=X
\]
 We define the pay-off 

\begin{equation}
J_{X,t}(v(.))=\int_{t}^{T}\mathbb{E}f(X(s),\mathcal{L}_{X(s)},v(X(s),s))ds+\mathbb{E}h(X(T),\mathcal{L}_{X(T)})\label{eq:9-7}
\end{equation}
and the value function 

\[
V(X,t)=\inf_{v(.)}\, J_{X,t}(v(.))
\]
 We claim that the value function is the solution of (\ref{eq:9-5}).
We proceed formally as above . We have to evaluate $V(X(t+\epsilon),t+\epsilon).$
But 

\[
X(t+\epsilon)\sim X+\epsilon g(X,\mathcal{L}_{X},v(X))+\sigma(w(t+\epsilon)-w(t))
\]
 Therefore 

\begin{align*}
V(X(t+\epsilon),t+\epsilon) & \sim\epsilon\frac{\partial V}{\partial t}(X,t)+\epsilon(DV(X,t),g(X,\mathcal{L}_{X},v))+\\
+ & V(X+\sigma(w(t+\epsilon)-w(t)),t)
\end{align*}
 and 

\[
V(X+\sigma(w(t+\epsilon)-w(t)),t)\sim V(X,t)+\frac{\sigma^{2}}{2}(D^{2}V(X,t)N,N)
\]
 We can then obtain (\ref{eq:9-5}) easily

\subsection{INTERPRETATION OF THE SECOND STOCHASTIC BELLMAN EQUATION}

By second stochastic Bellman equation we mean 

\begin{equation}
\frac{\partial V}{\partial t}+\frac{\sigma^{2}}{2}D^{2}V(X)(N,N)+\frac{1}{2}\beta^{2}\sum_{k=1}^{n}D^{2}V(X)(e_{k},e_{k})+\mathbb{E}H(X,\mathcal{L}_{X},DV(X))=0\label{eq:9-8}
\end{equation}
\[
V(X,T)=\mathbb{E}h(X,\mathcal{L}_{X})
\]
 We consider on $\Omega,\mathcal{A},$$P$ the variable $X,$ the
standard Wiener process $w(t)$ as above independent of $X$ and a
new standard Wiener process $b(t)$ with values in $\mathbb{R}^{n}$, independent
of $X$ and of $w(t).$ We set $\mathcal{B}^{t}=\sigma(b(s),s\leq t).$
A feedback control $v(X,s)$ on $\mathcal{H}$ is now a random field
on $\mathcal{H}$, such that, for $X$ fixed, it is a stochastic process
adapted to the filtration $\mathcal{B}^{s}$. We define the SDE on
$\mathcal{H}$ 

\begin{equation}
dX=g(X(s),\mathcal{L}_{X(s)},v(X(s),s))ds+\sigma dw(s)+\beta db(s)\label{eq:9-9}
\end{equation}
\[
X(t)=X
\]
 However, now $\mathcal{L}_{X(s)}$ does not refer to the probability
law of $X(s)$ but to the conditional probability given $\mathcal{B}^{s}$
. To figure out what it is, we consider the process $Y(s)$ defined
by the equation 

\begin{align}
dY & =g(Y(s)+\beta b(s),\mathcal{L}_{Y(s)+\beta b(s)},v(Y(s)+\beta b(s)))ds+\sigma dw(s)\label{eq:9-10}\\
Y(t) & =X-\beta b(t)\nonumber 
\end{align}

in which $b(s)$ must be considered as a given continuous function,
not a stochastic process. The only randomness in the model (\ref{eq:9-10})
is $X$ and $w(s).$ The equation (\ref{eq:9-10}) is similar to (\ref{eq:9-6})
since $b(s)$ is not stochastic. The conditional law of $X(s)$ given
$\mathcal{B}^{s}$ is $\mathcal{L}_{Y(s)+\beta b(s)}$ , so we can
write (\ref{eq:9-9}) as follows 

\begin{align}
dX & =g(X(s),\mathcal{L}_{Y(s)+\beta b(s)},v(X(s),s))ds+\sigma dw(s)+\beta db(s)\label{eq:9-11}\\
X(t) & =X\nonumber 
\end{align}

The payoff corresponding to the feedback $v(X(s),s)$ is now defined
by 

\begin{align}
J_{X,t}(v(.)) & =\mathbb{E}[\int_{t}^{T}f(X(s),\mathcal{L}_{Y(s)+\beta b(s)},v(X(s),s))ds+\label{eq:9-12}\\
+ & h(X(T),\mathcal{L}_{Y(T)+\beta b(T)})|\mathcal{B}^{t},X(t)=X]\nonumber 
\end{align}

The conditioning expesses the fact that we have access to the filtration
$\mathcal{B}^{t}$ . We define the value function 

\begin{equation}
V(X,t)=\inf_{v(.)}\, J_{X,t}(v(.))\label{eq:9-13}
\end{equation}

Normally $V(X,t)$ for fixed $X$ is random and $\mathcal{B}^{t}$
measurable. In fact, it will be independent of $\mathcal{B}^{t}$
and the solution of (\ref{eq:9-8}). This is due to the independence
of $X$ and $\mathcal{B}^{t}.$ We first write

\[
V(X,t)\sim\inf_{v(.)}\, \mathbb{E}\left[\epsilon f(X,\mathcal{L}_{X},v(X))\right.+
\]

\[
\left.\mathbb{E}\left(\int_{t+\epsilon}^{T}f(X(s),\mathcal{L}_{Y(s)+\beta b(s)},v(X(s),s))ds+h(X(T),\mathcal{L}_{Y(T)+\beta b(T)})|X(t+\epsilon),\mathcal{B}^{t+\epsilon},X(t)=X\right)|\mathcal{B}^{t},X(t)=X\right]
\]
 From the optimality principle, it follows 

\begin{equation}
V(X,t)\sim\inf_{v(.)}\, \mathbb{E}[\epsilon f(X,\mathcal{L}_{X},v(X))+\mathbb{E}(V(X(t+\epsilon),t+\epsilon)|\mathcal{B}^{t+\epsilon},X(t)=X)|\mathcal{B}^{t},X(t)=X]\label{eq:9-14}
\end{equation}
 Since 

\[
X(t+\epsilon)\sim X+\epsilon g(X,\mathcal{L}_{X},v(X))+\sigma(w(t+\epsilon)-w(t))+\beta(b(t+\epsilon)-b(t))
\]
 when we condition with respect to $\mathcal{B}^{t+\epsilon}$, we
must consider $b(t+\epsilon)-b(t)$ as fixed. Therefore we can check
easily that 

\[
\mathbb{E}(V(X(t+\epsilon),t+\epsilon)|\mathcal{B}^{t+\epsilon},X(t)=X)\sim V(X,t)+\epsilon\frac{\partial V}{\partial t}(X,t)+\epsilon(DV(X),g(X,\mathcal{L}_{X},v(X)))
\]

\[
+\epsilon\frac{\sigma^{2}}{2}(D^{2}V(X)N,N)+\beta\sum_{k=1}^{n}(DV(X),e_{k})(b_{k}(t+\epsilon)-b_{k}(t))+
\]
\[
+\frac{\beta^{2}}{2}\sum_{k,l=1}^{n}(D^{2}V(X)e_{k},e_{l})(b_{k}(t+\epsilon)-b_{k}(t))(b_{l}(t+\epsilon)-b_{l}(t))
\]
 in which $b_{k}(t)$ represent the coordinates of $b(t).$ Plugging
this formula in (\ref{eq:9-14}) and completing calculations we obtain
(\ref{eq:9-8}).

\section{THE MAXIMUM PRINCIPLE}

\subsection{PRELIMINARIES}

It is natural to consider the maximum principle in the context of
the new control problems considered in the previous section. We shall
need to change slightly the notation. We set $M=L^{2}(\Omega,\mathcal{A},P),$
hence $\mathcal{H}=M^{n}$. We should see $M$ as replacing the space
$\mathbb{R}$ of real numbers. We shall need $M^{d}$ for the space of controls.
A control is a function $V(s)$ with values in $M^{d}.$ When there
is no Wiener process $w(s)$ or $b(s)$, there is no aspect of adaptation.
So $V(s)$ is a ``deterministic'' function. Of course, this is paradoxical
since the values of $V(s)$ are random variables, but mathematically
the values are in a Hilbert space $M^{d}.$ When we introduce $w(t)$
and $b(t)$ we will say that it is a stochastic process.

\subsection{DETERMINISTIC MAXIMUM PRINCIPLE}

We consider the differential equation 

\begin{align}
\frac{dX}{dt} & =g(X(t),\mathcal{L}_{X(t)},v(t))\label{eq:10-1}\\
X(0) & =X\nonumber 
\end{align}

in which the control $v(t)$ takes values in $M^{d}.$The solution
$X(t)$ takes values in $M^{n}$. The particularity of the space $M^{n}$
is that we may consider the random variable together with its probability
distribution .

We then define the payoff functional 

\begin{equation}
J(v(.))=\int_{0}^{T}\mathbb{E}f(X(t),\mathcal{L}_{X(t)},v(t))dt+\mathbb{E}h(X(T),\mathcal{L}_{X(T)})\label{eq:10-2}
\end{equation}

The control $v(t)$ can be defined by a feedback $v(t)=v(X(t),\mathcal{L}_{X(t)},t).$
We recall the Dynamic Programming approach to this problem . We have 

\begin{equation}
\inf_{v(.)}\, J(v(.))=V(X,0)\label{eq:10-3}
\end{equation}

in which the value function $V(X,t)$ is the solution of the first
order equation 

\begin{align}
\frac{\partial V}{\partial t}+\mathbb{E}H(X,\mathcal{L}_{X},DV(X)) & =0\label{eq:10-4}\\
V(X,T)=\mathbb{E}\, h(X,\mathcal{L}_{X})\nonumber 
\end{align}

in which we recall the notation 

\[
H(x,m,q)=\inf_{v}(f(x,m,v)+q.g(x,m,v))
\]

and we have denoted by $\hat{v}(x,m,q)$ the point of minimum. Next,
consider for $Q$ in $\mathcal{H}=M^{n},$the new Hamiltonian 

\begin{equation}
\tilde{H}(X,\mathcal{L}_{X},Q)=\inf_{v\in M^{d}}[\mathbb{E}f(X,\mathcal{L}_{X},v)+(Q,g(X,\mathcal{L}_{X},v))]\label{eq:10-5}
\end{equation}

in view of the definition of the scalar product the infimum is attained
at $\hat{v}(X,\mathcal{L}_{X},Q)$. Therefore 

\[
\tilde{H}(X,\mathcal{L}_{X},Q)=\mathbb{E}H(X,\mathcal{L}_{X},Q)
\]

We next define the Lagrangian 

\begin{equation}
L(X,\mathcal{L}_{X},v,Q)=\mathbb{E}f(X,\mathcal{L}_{X},v)+(Q,g(X,\mathcal{L}_{X},v))\label{eq:10-6}
\end{equation}

in which $X,Q\in M^{n}$ and $v\in M^{d}$. There is a slight abuse
of notation for $v$ which is an element of $\mathbb{R}^{d}$ in the usual
defintion of the Hamiltonian and an element of $M^{d}$ in the definition
of the Lagrangian. We have 

\begin{equation}
\tilde{H}(X,\mathcal{L}_{X},Q)=\inf_{v\in M^{d}}L(X,\mathcal{L}_{X},v,Q)\label{eq:10-7}
\end{equation}

and the infimum is attained at $\hat{v}(X,\mathcal{L}_{X},Q).$ Bellman
equation (\ref{eq:10-4}) reads 

\begin{align}
\frac{\partial V}{\partial t}+\tilde{H}(X,\mathcal{L}_{X},DV(X)) & =0\label{eq:10-8}\\
V(X,T)=\mathbb{E}\, h(X,\mathcal{L}_{X})\nonumber 
\end{align}

We note that $V(X,t)=V(\mathcal{L}_{X},t)$ so $DV(X,t)=D_{x}\dfrac{\partial V(\mathcal{L}_{X},t)}{\partial m}(X).$
Therefore 

\[
\hat{v}(X,\mathcal{L}_{X},DV(X,t))=\hat{v}(X,\mathcal{L}_{X},D_{x}\dfrac{\partial V(\mathcal{L}_{X},t)}{\partial m}(X))
\]
 which is a feedback in $X$ as defined above. We set also 

\begin{align}
F(X,\mathcal{L}_{X},Q) & =f(X,\mathcal{L}_{X},\hat{v}(X,\mathcal{L}_{X},Q))\label{eq:10-9}\\
G(X,\mathcal{L}_{X},Q) & =g(X,\mathcal{L}_{X},\hat{v}(X,\mathcal{L}_{X},Q))\nonumber 
\end{align}

and 

\[
H(X,\mathcal{L}_{X},Q)=F(X,\mathcal{L}_{X},Q)+Q.G(X,\mathcal{L}_{X},Q)
\]

We shall need $D_{X}\tilde{H}(X,\mathcal{L}_{X},Q).$ By that we mean
that $Q$ is a fixed element of $M^{n}$ and we take the gradient
with respect to $X.$ Since $\tilde{H}(X,\mathcal{L}_{X},Q)=\mathbb{E}H(X,\mathcal{L}_{X},Q)$
we can compute the gradient according to formula (\ref{eq:1-32})
which means 

\begin{equation}
D_{X}\tilde{H}(X,\mathcal{L}_{X},Q)=D_{x}H(X,\mathcal{L}_{X},Q)+D_{x}\mathbb{E}_{YQ}\frac{\partial H}{\partial m}(Y,\mathcal{L}_{X},Q)(X)\label{eq:10-10}
\end{equation}

in which $Y$ is a copy of $X$ . We have also 

\begin{align}
D_{X}\tilde{H}(X,\mathcal{L}_{X},Q) & =D_{x}F(X,\mathcal{L}_{X},Q)+(D_{x}G)^{*}(X,\mathcal{L}_{X},Q)Q+\label{eq:10-11}\\
+D_{x}\mathbb{E}_{YQ}\frac{\partial F(Y,\mathcal{L}_{X},Q)}{\partial m}(X) & +\mathbb{E}_{YQ}(D_{x}\frac{\partial G}{\partial m})^{*}(Y,\mathcal{L}_{X},Q)(X)Q\nonumber 
\end{align}
 Set $\mathcal{U}(X,t)=DV(X,t),$ then Bellman equation reads 

\begin{equation}
\frac{\partial V}{\partial t}+\tilde{H}(X,\mathcal{L}_{X},\mathcal{U}(X,t))=0\label{eq:10-12}
\end{equation}

Differentiating in $X$ we obtain the Master equation 

\begin{equation}
\frac{\partial\mathcal{U}}{\partial t}+D_{X}\tilde{H}(X,\mathcal{L}_{X},\mathcal{U})+D\mathcal{U}(X)G(X,\mathcal{L}_{X},\mathcal{U})=0\label{eq:10-13}
\end{equation}

Consider the optimal control problem (\ref{eq:10-1}), (\ref{eq:10-2}),
the optimal control $u(t)$ is given by 

\begin{equation}
u(t)=\hat{v}(X(t),\mathcal{L}_{X(t)},\mathcal{U}(X(t),t))\label{eq:10-14}
\end{equation}
 and using the same notation $X(t)$ for the optimal state we have

\begin{equation}
g(X(t),\mathcal{L}_{X(t)},u(t))=G(X(t),\mathcal{L}_{X(t)},\mathcal{U}(X(t),t))\label{eq:10-15}
\end{equation}
\begin{equation}
L(X(t),\mathcal{L}_{X(t)},u(t),\mathcal{U}(X(t),t))=\tilde{H}(X(t),\mathcal{L}_{X(t)},\mathcal{U}(X(t),t))\label{eq:10-16}
\end{equation}

Define now $Z(t)=\mathcal{U}(X(t),t),$ the co state, we have 

\[
\frac{dZ(t)}{dt}=\frac{\partial\mathcal{U}(X(t),t)}{\partial t}+D\mathcal{U}(X(t),t)G(X(t),\mathcal{L}_{X(t)},\mathcal{U}(X(t),t))
\]
 so, collecting results from (\ref{eq:10-13}) to (\ref{eq:10-16})
we obtain the following system 

\[
\frac{dX}{dt}=g(X(t),\mathcal{L}_{X(t)},u(t))
\]

\[
-\frac{dZ(t)}{dt}=D_{X}L(X(t),\mathcal{L}_{X(t)},u(t),Z(t))
\]

\begin{equation}
X(0)=X\label{eq:10-17}
\end{equation}

\[
Z(T)=D_{X}h(X(T),\mathcal{L}_{X(T)})
\]
\begin{equation}
u(t)\:\text{minimizes }L(X(t),\mathcal{L}_{X(t)},v,Z(t))\:\text{in}\: v\in M^{d}\label{eq:10-18}
\end{equation}

where the Lagrangian is defined in (\ref{eq:10-6}). For completion
, we express 

\[
D_{X}L(X(t),\mathcal{L}_{X(t)},u(t),Z(t))=D_{x}f(X(t),\mathcal{L}_{X(t)},u(t))+\mathbb{E}_{Y(t)u(t)}D_{x}\frac{\partial}{\partial m}f(Y(t),\mathcal{L}_{X(t)},u(t))(X(t))+
\]
\begin{equation}
+(D_{x}g)^{*}(X(t),\mathcal{L}_{X(t)},u(t))Z(t)+\mathbb{E}_{Y(t)Z(t)u(t)}(D_{x}\frac{\partial g}{\partial m})^{*}(Y(t),\mathcal{L}_{X(t)},u(t))(X(t))Z(t)\label{eq:10-19}
\end{equation}
 in which, as usual $Y(t)$ is a copy of $X(t).$ For the sake of
verification, let us compute $(D_{X}L(X(t),\mathcal{L}_{X(t)},u(t),Z(t)),\tilde{X}(t))$in
two ways. First from the definition of the Lagrangian we can write 

\[
(D_{X}L(X(t),\mathcal{L}_{X(t)},u(t),Z(t)),\tilde{X}(t))=(D_{X}\mathbb{E}f(X(t),\mathcal{L}_{X(t)},u(t)),\tilde{X}(t))+
\]

\[
+(D_{X}(g(X(t),\mathcal{L}_{X(t)},u(t)),Z(t)),\tilde{X}(t))
\]
 From formula (\ref{eq:1-32}) we see that $D_{X}\mathbb{E}f(X(t),\mathcal{L}_{X(t)},u(t))$
is indeed equal to the term in $f$ in the right-hand side of (\ref{eq:10-19}).
Next 

\[
(D_{X}(g(X(t),\mathcal{L}_{X(t)},u(t)),Z(t)),\tilde{X}(t))=(D_{X}g(X(t),\mathcal{L}_{X(t)},u(t))\tilde{X}(t),Z(t))
\]
 From formula (\ref{eq:1-35}) we can write 

\begin{align*}
D_{X}g(X(t),\mathcal{L}_{X(t)},u(t))\tilde{X}(t) & =D_{x}g(X(t),\mathcal{L}_{X(t)},u(t))\tilde{X}(t)+
\end{align*}

\[
+\mathbb{E}_{Y(t)\tilde{X}(t)}D_{y}\frac{\partial}{\partial m}g(X(t),\mathcal{L}_{X(t)},u(t))(Y(t))\tilde{X}(t)
\]
 therefore 

\[
(D_{X}g(X(t),\mathcal{L}_{X(t)},u(t))\tilde{X}(t),Z(t))=\mathbb{E}\sum_{ij}Z_{i}(t)\frac{\partial g_{i}}{\partial x_{j}}(X(t),\mathcal{L}_{X(t)},u(t))\tilde{X_{j}}(t)+
\]

\[
+\mathbb{E}_{Z(t)X(t)u(t)}\sum_{i}Z_{i}(t)\mathbb{E}_{Y(t)\tilde{X}(t)}\dfrac{\partial}{\partial y_{j}}\frac{\partial}{\partial m}g_{i}(X(t),\mathcal{L}_{X(t)},u(t))(Y(t))\tilde{X_{j}}(t)=I+II
\]
 Consider the second term, which is the only point to check.We can
exchange the names of $X(t)$ and $Y(t)$ which are identical, hence 

\[
II=\mathbb{E}_{X(t)\tilde{X}(t)}\mathbb{E}_{Z(t)Y(t)u(t)}\sum_{i}Z_{i}(t)\dfrac{\partial}{\partial x_{j}}\frac{\partial}{\partial m}g_{i}(Y(t),\mathcal{L}_{X(t)},u(t))(X(t))\tilde{X_{j}}(t)
\]

Now, referring to (\ref{eq:10-19}) and testing the formula with $\tilde{X}(t),$we
obtain a formula which is identical to $I+II.$ So we can state

\begin{equation}
(D_{X}L(X(t),\mathcal{L}_{X(t)},u(t),Z(t)),\tilde{X}(t))=(D_{X}\mathbb{E}f(X(t),\mathcal{L}_{X(t)},u(t)),\tilde{X}(t))+\label{eq:10-20}
\end{equation}

\[
+(D_{X}g(X(t),\mathcal{L}_{X(t)},u(t))\tilde{X}(t),Z(t))
\]
 
\begin{xca}
\label{exer5-1}We can derive the system (\ref{eq:10-18}) directly,
not as a consequence of the Bellman equation. 
\end{xca}

\subsection{STOCHASTIC MAXIMUM PRINCIPLE}

We introduce the Wiener process $w(t)$ with values in $\mathbb{R}^{n},$ independent
of the initial condition $X$ and the filtration $\mathcal{F}^{t}$
generated by $X$ and $w(s),s\leq t.$ A control $v(t)$ is a stochastic
process adapted to $\mathcal{F}^{t}$ with values in $M^{d}.$ We
then define the state of the system by the SDE in $M^{d}$ 

\begin{align}
dX & =g(X(t),\mathcal{L}_{X(t)},v(t))dt+\sigma dw(t)\label{eq:10-21}\\
X(0) & =X\nonumber 
\end{align}

and the pay off is defined by 

\begin{equation}
J(v(.))=\int_{0}^{T}\mathbb{E}f(X(t),\mathcal{L}_{X(t)},v(t))dt+\mathbb{E}h(X(T),\mathcal{L}_{X(T)})\label{eq:10-22}
\end{equation}

We know that 

\begin{equation}
\inf_{v(.)}J(v(.))=V(X,0)\label{eq:10-23}
\end{equation}

in which $V(X,t)$ is the solution of the Bellman equation 

\begin{align}
\frac{\partial V}{\partial t}+\frac{\sigma^{2}}{2}D^{2}V(X)(N,N)+\tilde{H}(X,\mathcal{L}_{X},DV(X)) & =0\label{eq:10-24}\\
V(X,T)=\mathbb{E}\, h(X,\mathcal{L}_{X})\nonumber 
\end{align}

The optimal control $u(t)$ is obtained by minimising the Lagrangian
( see \ref{eq:10-6}) $L(X(t),\mathcal{L}_{X(t)},v,\mathcal{U}(X(t),t))$
in $v\in M^{d},$ with $\mathcal{U}(X(t),t)=DV(X(t),t).$ The function
$\mathcal{U}(X,t)$ is the solution of the Master equation

\begin{equation}
\frac{\partial\mathcal{U}}{\partial t}+\frac{\sigma^{2}}{2}(D^{2}\mathcal{U}(X)N,N)+D_{X}\tilde{H}(X,\mathcal{L}_{X},\mathcal{U})+D\mathcal{U}(X)G(X,\mathcal{L}_{X},\mathcal{U})=0\label{eq:10-25}
\end{equation}

Define $Z(t)=\mathcal{U}(X(t),t)$ , where $X(t)$ is the optimal
state ( we use the same notation), then, by Ito's calculus is the
space $M^{n}$ we have 

\[
dZ(t)=D\mathcal{U}(X(t),t)g(X(t),\mathcal{L}_{X(t)},u(t))dt+\sigma D\mathcal{U}(X(t),t)dw(t)+\frac{\sigma^{2}}{2}(D^{2}\mathcal{U}(X(t))N,N)dt+\frac{\partial\mathcal{U}}{\partial t}(X(t),t)dt
\]
 therefore from (\ref{eq:10-25}) it follows, setting $K(t)=$$\sigma D\mathcal{U}(X(t),t),$ 

\begin{align}
-dZ(t) & =D_{X}L(X(t),\mathcal{L}_{X(t)},u(t),Z(t))dt-K(t)dw(t)\label{eq:10-26}\\
Z(T) & =D_{X}h(X(T),\mathcal{L}_{X(T)})\nonumber 
\end{align}

We have also the state equation

\begin{align}
dX(t) & =g(X(t),\mathcal{L}_{X(t)},u(t))dt+\sigma dw(t)\label{eq:10-27}\\
X(0) & =X\nonumber 
\end{align}
 and the optimality condition for $u(t)$ is expressed by 

\begin{equation}
u(t)\:\text{minimizes }L(X(t),\mathcal{L}_{X(t)},v,Z(t))\:\text{in}\: v\in M^{d}\label{eq:10-28}
\end{equation}

The triple (\ref{eq:10-26}),(\ref{eq:10-27}),(\ref{eq:10-28}) expresses
the stochastic maximum principle. The adjoint equation (\ref{eq:10-26})
is a BSDE in the space $M^{d}.$

\end{document}